\documentclass[authoryear, 11pt]{imsart}
\RequirePackage[OT1]{fontenc}
\RequirePackage{amsthm,amsmath,natbib}
\RequirePackage{amsthm,amsmath}
\RequirePackage{hypernat}
\startlocaldefs
\numberwithin{equation}{section}
\theoremstyle{plain}

\endlocaldefs
\usepackage{amsfonts,amsmath,latexsym,amssymb}
\usepackage{fullpage}

\newtheorem{lemma}{Lemma}
\newtheorem{theorem}{Theorem}
\newtheorem{proposition}{Proposition}

\newtheorem{assumption}{Assumption}
\newtheorem{remark}{Remark}

\newtheorem{corollary}{Corollary}

\renewcommand{\kappa}{\varkappa}

\newcommand{\rd}{{\rm d}}
\newcommand{\e}{\varepsilon}

\newcommand{\cA}{{\cal A}}

\newcommand{\cK}{{\cal K}}
\newcommand{\cL}{{\cal L}}

\newcommand{\cQ}{{\cal Q}}

\newcommand{\cS}{{\cal S}}

\newcommand{\cX}{{\cal X}}

\newcommand{\bB}{\mathbb B}
\newcommand{\bC}{\mathbb C}

\newcommand{\bE}{\mathbb E}

\newcommand{\bH}{\mathbb H}

\newcommand{\bL}{{\mathbb L}}

\newcommand{\bN}{{\mathbb N}}

\newcommand{\bP}{{\mathbb P}}

\newcommand{\bR}{{\mathbb R}}
\newcommand{\bS}{{\mathbb S}}
\newcommand{\bT}{{\mathbb T}}

\newcommand{\mA}{\mathfrak{A}}
\newcommand{\mB}{\mathfrak{B}}
\newcommand{\mE}{\mathfrak{E}}

\newcommand{\mS}{\mathfrak{S}}

\newcommand{\mh}{\mathfrak{h}}

\newcommand{\mT}{\mathfrak{T}}

\newcommand{\mH}{\mathfrak{H}}

\newcommand{\epr}{\hfill\hbox{\hskip 4pt
                \vrule width 5pt height 6pt depth 1.5pt}\vspace{0.5cm}\par}
\def    \qed    {\hfill\hbox{\hskip 4pt
                \vrule width 3pt height 1pt depth 1pt
               \hbox{\vrule width 2pt height 3.5pt depth 1pt}}}

\begin{document}
\bibliographystyle{plainnat}
\begin{frontmatter}
\title{Upper functions for $L_p$-norm of gaussian random fields}
\runtitle{Upper functions}

\begin{aug}
\author{\fnms{Oleg} \snm{Lepski}
\ead[label=e2]{lepski@cmi.univ-mrs.fr}}
\runauthor{ O. Lepski}

\affiliation{Aix--Marseille  Universit\'e}


\address{Laboratoire d'Analyse, Topologie, Probabilit\'es\\
  Aix-Marseille Universit\'e \\
 39, rue F. Joliot-Curie \\
13453 Marseille, France\\
\printead{e2}\\ }
\end{aug}

\begin{abstract}
In this paper we
are interested in finding upper functions for a collection of  random variables
$\big\{\big\|\xi_{\vec{h}}\big\|_p, \vec{h}\in\mathrm{H}\big\}, 1\leq p<\infty$.
Here
$\xi_{\vec{h}}(x), x\in(-b,b)^d, d\geq 1$ is a kernel-type gaussian random field and
$\|\cdot\|_p$ stands for $\bL_p$-norm on $(-b,b)^d$. The set $\mathrm{H}$ consists of $d$-variate vector-functions
defined on $(-b,b)^d$ and taking values in some countable net in $\bR^d_+$.
We seek a non-random family $\left\{\Psi_\alpha\big(\vec{h}\big),\;\;\vec{h}\in\mathrm{H}\right\}$ such that
$
\bE\big\{\sup_{\vec{h}\in\mathrm{H}}\big[\big\|\xi_{\vec{h}}\big\|_p-\Psi_\alpha\big(\vec{h}\big)\big]_+\big\}^q\leq \alpha^q,\; q\geq 1,
$
where $\alpha>0$ is prescribed level.

\end{abstract}


\maketitle

\end{frontmatter}

\section{Introduction}
\label{sec:introduction}
Let  $\bR^d,\;d\geq 1,$ be equipped with  Borel $\sigma$-algebra $\mB(\bR^d)$  and Lebesgue measure $\nu_d$. Put $\widetilde{\mB}(\bR^d)=\left\{B\in\mB(\bR^d):\;\;\nu_d(B)<\infty\right\}$ and let $\big(W_B,\; B\in\widetilde{\mB}(\bR^d)\big)$ be the white noise with intensity $\nu_d$.
Throughout of the paper we will use the following notations. For any $u,v\in\bR^d$ the operations and relations $u/v$, $uv$, $u\vee v$,$u\wedge v$,
$u<v$, $au, a\in\bR,$ are understood in coordinate-wise sense and $|u|$ stands for euclidian norm of $u$.
All integrals are taken over $\bR^d$
unless the domain of integration is specified explicitly. For any real $a$ its  positive part is denoted by  $(a)_+$ and $\lfloor a\rfloor$ is used for its integer part. For any $\mathbf{n}=(n_1,\ldots,n_d)\in\bN^d,\;d\geq 1,$  $|\mathbf{n}|$ stands for $\sum_{j=1}^dn_i$.

\paragraph{Collection of random variables}

Let    $0<\mh\leq e^{-2}$ be fixed number and put $\mH=\{\mh_s, s\in\bN\}$, where $\mh_s=e^{-s}\mh$.
Denote by  $\mathfrak{S}(\mh)$ the set of all measurable functions defined on $(-b,b)^d,\;b\in (0,\infty),$ and  taking values in $\mH$
and define
$$
\mathfrak{S}_d(\mh)=\Big\{\vec{h}:(-b,b)^{d}\to\mH^d:\quad \vec{h}(x)=\big(h_1(x),\ldots,h_d(x)\big),\;x\in(-b,b)^d,\;\;h_i\in\mS(\mH),\;i=\overline{1,d}\Big\}.
$$
Let $K:\bR^d\to\bR$ be fixed. With any
$\vec{h}\in\mathfrak{S}_d(\mh)$ we associate the function
$$
K_{\vec{h}}(t,x)=V^{-1}_{\vec{h}}(x)K\bigg(\frac{t-x}{\vec{h}(x)}\bigg),\;\; t\in\bR^d,\; x\in(-b,b)^d,
$$
where $V_{\vec{h}}(x)=\prod_{i=1}^dh_i(x)$. Following the terminology used in the mathematical statistics we call the function $K$ kernel and the vector-function  $\vec{h}$ multi-bandwidth. Moreover, if all coordinates of $\vec{h}$ are the same  we will say that corresponding collection
 is {\it isotropic}. Otherwise it is called {\it anisotropic}.

Let $\mathrm{H}$ be a given subset of $\mathfrak{S}_d(\mh)$ and  consider the family
$$
\left\{\xi_{\vec{h}}(x)=\int K_{\vec{h}}(t,x)W(\rd t),\;\;\vec{h}\in\mathrm{H},\;x\in(-b,b)^d\right\}.
$$
We note that $\xi_{\vec{h}}$ is centered  gaussian random field on $(-b,b)^d$ with the covariance function
$$
V^{-1}_{\vec{h}}(x)V^{-1}_{\vec{h}}(y)\int K\bigg(\frac{t-x}{\vec{h}(x)}\bigg)K\bigg(\frac{t-y}{\vec{h}(y)}\bigg)\nu_d(\rd t),\quad x,y\in(-b,b)^d.
$$
Throughout the paper $(\xi_{\vec{h}}, \vec{h}\in\mathrm{H})$ is supposed to be defined on the probability space $(\mathfrak{X},\mA,\bP)$ and furthermore $\bE$ denotes the expectation with respect to $\bP$.

\paragraph{Objectives} Our goal is to find an upper function for the following collection of random variables
$$
\Lambda_p\left(\mathrm{H}\right)=\left\{\big\|\xi_{\vec{h}}\big\|_p,\;\;\vec{h}\in\mathrm{H}\right\},\quad 1\leq p<\infty,
$$
 where  $\|\cdot\|_p$ stands for $\bL_p$-norm on $(-b,b)^d$.

 More precisely we seek for a {\it non-random collection} $\left\{\Psi_\alpha\big(\vec{h}\big),\;\;\vec{h}\in\mathrm{H}\right\}$ such that
\begin{equation}
\label{eq:objective}
\bE\bigg\{\sup_{\vec{h}\in\mathrm{H}}\Big[\big\|\xi_{\vec{h}}\big\|_p-\Psi_\alpha\big(\vec{h}\big)\Big]_+\bigg\}^q\leq \alpha^q,\quad q\geq 1,
\end{equation}
where $\alpha >0$ is a prescribed level.

It is worth mentioning that
uniform probability and moment bounds for $\left[\sup_{\theta\in\Theta}\Upsilon(\chi_\theta)\right]$ in the case where
$\chi_\theta$ is
 empirical or gaussian  process and $\Upsilon$  is a positive functional are a subject of vast literature, see, e.g.,
\cite{alexander},  \cite{Talagrand,Tal-book}, \cite{Lif}, \cite{wellner}, \cite{van-de-Geer},
\cite{massart}, \cite{bousquet},
\cite{gine-kolt} among many others.
Such bounds play an important role in establishing
the laws of iterative logarithm and central limit
theorems [see, e.g., \cite{alexander} and
\cite{gine-zinn}].
However much less attention was paid to finding of  upper functions.
Some asymptotical results can be found in \cite{Kalina}, \cite{Qua&Vatan}, \cite{Bobkov}, \cite{Shiryaev} and references therein.
The inequalities similar to (\ref{eq:objective}) was obtained \cite{Ostrovskii}, \cite{GL3} and \cite{lepski-a,lepski-b,lepski-c}.

The evaluation of upper functions have become the important technical tool in different areas of mathematical  statistics in particular in adaptive estimation. Indeed,
 almost all known constructions of adaptive estimators e.g. \cite{Barron-Birge}, \cite{golubev},   \cite{GL2,GL11}
 involve the computation of upper functions for stochastic objects of different kinds. We provide below with explicit expression of the functional
$\Psi_\alpha$ that allows, in particular,  to use our results for constructing data-driven  procedures in multivariate function estimation.

The upper functions for $\bL_p$-norm of "kernel-type" empirical and gaussian processes  was studied in recent papers \cite{GL3} and \cite{lepski-a}.
However the results obtained there allow to study only bandwidth's collection consisted of constant functions, see discussions after Theorems \ref{th:deviation-l_p-gauss-general}--\ref{th:deviation-l_p-gauss-isotropic} below. To the best of our knowledge the problem of constructing  upper functions
for the collection parameterized by bandwidths being multivariate functions was not studied in the literature.

\paragraph{Assumptions imposed on the kernel $K$}
Throughout the paper we will consider the collections $\Lambda(\mathrm{H})$ with $K$ satisfying one of  assumptions indicated below.
Let $a\geq 1$ and $L>0$ be fixed.
\begin{assumption}
\label{ass:kernel-new-new}
  $\text{supp}(K)\subset[-a,a]^d$ and
\begin{eqnarray*}
|K(s)-K(t)|\leq L|s-t|,\;\;\forall s,t\in\bR^d;
\end{eqnarray*}
\end{assumption}

\begin{assumption}
\label{ass:kernel-new}
  $\text{supp}(K)\subset[-a,a]^d$ and for any $\mathbf{n}\in\bN$ such that $|\mathbf{n}|\leq \lfloor d/2\rfloor+1$
\begin{eqnarray*}
|D^{\mathbf{n}}K(s)-D^{\mathbf{n}}K(t)|\leq L|s-t|,\;\;\forall s,t\in\bR^d,
\end{eqnarray*}
where   $D^{\mathbf{n}}=\frac{\partial^{|\mathbf{n}|}}{\partial y_1^{n_1}\cdots\partial y_k^{n_d}}.$
\end{assumption}

\begin{assumption}
\label{ass:kernel}
 There exists
  $\cK:\bR\to\bR$  such that $\text{supp}(K)\subset[-a,a]$ and
\begin{eqnarray*}
&(\mathbf{i})&\qquad|\cK(s)-\cK(t)|\leq L|s-t|,\;\;\forall s,t\in\bR;
\\
&(\mathbf{ii})&\qquad K(x)=\prod_{i=1}^d\cK(x_i),\;\;\forall x=(x_1,\ldots,x_d)\in\bR^d.
\end{eqnarray*}

\end{assumption}

\paragraph{Organization of the paper}

In Section~\ref{sec:main-results} we present three constructions of upper functions and proved for them
the inequality of type (\ref{eq:objective}), Theorems \ref{th:deviation-l_p-gauss-general}--\ref{th:deviation-l_p-gauss-isotropic}.
Moreover, in Subsection \ref{sec:examples} we discuss the example of the bandwidth collection satisfying  assumptions of Theorem \ref{th:deviation-l_p-gauss}. Section \ref{sec:proofs} contains
proofs of Theorems~\ref{th:deviation-l_p-gauss-general}--\ref{th:deviation-l_p-gauss-isotropic}; proofs of auxiliary results are relegated to Appendix.

\section{Main results}
\label{sec:main-results}

\subsection{Anisotropic case. First construction}
\label{sec:selection-rule-and-L_p-norm-oracle-inequality}

For any    $\vec{h}\in\mS_d(\mh)$ and any $0<\e\leq e^{-2}$ define
\begin{eqnarray*}
&&\psi_\e\big(\vec{h}\big)=C_1\Big\|\sqrt{\big|\ln{\big(\e V_{\vec{h}}\big)}\big|} V^{-\frac{1}{2}}_{\vec{h}}\Big\|_p,
\end{eqnarray*}
where
$
C_1=2(q\vee p)+2\sqrt{2d}\Big[\sqrt{\pi}+ \|K\|_2\Big(\sqrt{\big|\ln{\big(4bL\|K\|_2\big)}\big|} +1\Big)\Big].
$
\begin{theorem}
\label{th:deviation-l_p-gauss-general}
Let  $q\geq 1$, $p\geq 1$,  be fixed and let $\mathrm{H}$ be an arbitrary
countable subset of $\mS_d(\mh)$. Suppose also that Assumption \ref{ass:kernel-new-new} is fulfilled.
\begin{eqnarray*}
&&\bE\bigg\{\sup_{\vec{h}\in\mathrm{H}}
\Big[\big\|\xi_{\vec{h}}\big\|_p-\psi_\e(\vec{h})\Big]_+\bigg\}^{q}
\leq \big[C_3\e\big]^{q}, \quad\forall\mh,\e\in \big(0,e^{-2}\big),
\end{eqnarray*}
where
$
C_3=2^{\frac{d}{p}}\bigg[2\tilde{q}\int_0^\infty z^{\tilde{q}-1}\exp\bigg(-\frac{z^{\frac{2}{p}}}{8\|K\|^{2}_2}\bigg)\rd z\bigg]^{\frac{1}{p\tilde{q}}},
\quad \tilde{q}=(q/p)\vee 1.
$
\end{theorem}

\begin{remark}
\label{rem:after-th:deviation-l_p-gauss-general}
We consider only countable subsets of $\mS_d(\mh)$ in order not to discuss the measurability issue. Actually the statement of the theorem remains valid for any subset providing the measurability of  the corresponding  supremum. It explains why the upper function $\psi_\e$ as well as the constants $C_1$ and $C_3$
are independent of the choice of $\mathrm{H}$.
\end{remark}

The advantage of the result presented in Theorem \ref{th:deviation-l_p-gauss-general}  is that it is proved without any condition imposed on the set of bandwidths. However the natural question arising in this context is the presented bound sharp whatever the choice of $\mathrm{H}$? We will discuss this issue more in detail in the next section. Here we only say that the answer on the aforementioned question is negative. Indeed, let us suppose that
$h_i(x)=h_j(x), i,j=1,\ldots d$ (isotropic case) and additionally  $h_1(x)=h,\; h\in\mH$, for any $x\in (-b,b)^d$ .
In other words the bandwidths we consider are the constants.
In this case
$$
\psi_\e(\vec{h})=C_1\sqrt{|\ln(\e)|+d|\ln(h)|}h^{-\frac{d}{2}}.
$$
However, the upper function found  in \cite{lepski-a}, Theorem 1,  is given by
$$
\Psi(h)=Ch^{-\frac{d}{2}},\quad \forall p\geq 2.
$$
The level provided  by this upper function is also proportional  $\e^q$,  under assumption $\mh\leq \big(c|\ln(\e)|\big)^{\frac{p}{2}}$.  Here $C$ and $c$ are absolute constants. Moreover, as it was mentioned in \cite{lepski-a} if in the considered example  $\mathrm{H}$ consists of a single element $h$ then
$$
\Big( \bE
\big\|\xi_{\vec{h}}\big\|^q_p\Big)^{\frac{1}{q}}\geq c_p h^{-\frac{d}{2}}.
$$
Thus we can assert that $\Psi(h)$ is sharp.

As we see $\psi_\e(\vec{h})\gg \Psi(h)$ for all $h\in\mathrm{H}$. The question  we address now is: can the upper function given in Theorem \ref{th:deviation-l_p-gauss-general} be improved when an arbitrary collection $\mathrm{H}$ is considered? Our conjecture that the answer is  negative in general but for  sets of bandwidths satisfying rather weak assumption presented below it is possible.

\subsection{Anisotropic case. Functional classes of bandwidths}
\label{sec:section-upper-functions}


Put for any $\vec{h}\in\mathfrak{S}_d(\mh)$ and any multi-index $\mathbf{s}=(s_1,\ldots,s_d)\in\bN^d$
$$
\Lambda_{\mathbf{s}}\big[\vec{h}\big]=\cap_{j=1}^{d}\Lambda_{s_j}\big[h_j\big],\qquad \Lambda_{s_j}\big[h_j\big]=\big\{x\in(-b,b)^d:\;\; h_j(x)=\mh_{s_j}\big\}.
$$
Let $\tau\in (0,1)$ and $\cL>0$ be given constants. Define
\begin{eqnarray*}
\bH_d(\tau,\cL)=\bigg\{\vec{h}\in\mathfrak{S}_d(\mh):\;\; \sum_{\mathbf{s}\in\cS_d}\nu_d^{\tau}\Big(\Lambda_{\mathbf{s}}\big[\vec{h}\big]\Big)\leq \cL\bigg\}.
\end{eqnarray*}
Put $\bN^*_p=\big\{[p]+1,[p]+2,\ldots\big\}$ and introduce for any $\mathcal{A}\geq \mh^{-\frac{d}{2}}$
$$
\mathbb{B}(\cA)=\bigcup_{r\in \bN^*_p}\mathbb{B}_r(\mathcal{A}),\qquad \mathbb{B}_r(\mathcal{A})=\bigg\{\vec{h}\in\mathfrak{S}_d(\mh):\;\;
\Big\| V^{-\frac{1}{2}}_{\vec{h}}\Big\|_{\frac{rp}{r-p}}\leq \cA\bigg\}.
$$
In this section we will be interested in finding an upper function  when $\mathrm{H}$ is an arbitrary subset
of $\bH_d(\tau,\cL,\cA):=\bH_d(\tau,\cL)\cap\mathbb{B}(\cA)$.

A  simple example of the subset of $\bH_d(\tau,\cL,\cA)$ is given in paragraph $\mathbf{2^0}$ below while a quite sophisticated construction is postponed to  Section \ref{sec:examples}.

\smallskip

The following notations related to the functional class $\mathbb{B}(\cA)$ will be exploited in the sequel. For any $\vec{h}\in\mathbb{B}(\cA)$ define
\begin{equation}
\label{eq:def-N_p(h)}
\bN^*_p\big(\vec{h},\cA\big)= \bN^*_p\cap \big[r_{\cA}(\vec{h}),\infty\big),
\qquad
r_{\cA}(\vec{h})=\inf\big\{r\in \bN^*_p:\;\; \vec{h}\in\mathbb{B}_r(\mathcal{A})\big\}.
\end{equation}
Obviously   $r_{\cA}\big(\vec{h}\big)<\infty$ for any $\vec{h}\in\mathbb{B}(\mathcal{A})$.


The following relation between parameters $\mh, \cA$ and $\tau$ is supposed to be held throughout of this section.
\begin{equation}
\label{eq0:relations-between-parameters}
d\ln\ln(\cA)\leq 2\sqrt{2(1-\tau)|\ln(\mh)|}-d\ln(4).
\end{equation}
For any    $\vec{h}\in\bB(\cA)$ define
\begin{eqnarray*}
&&\psi\big(\vec{h}\big)=
\inf_{r\in\bN^*_p(\vec{h})}C_2(r,\tau,\cL)\Big\| V^{-\frac{1}{2}}_{\vec{h}}\Big\|_{\frac{rp}{r-p}},
\end{eqnarray*}
where,  $\bN^*_p(\vec{h})$ is defined in (\ref{eq:def-N_p(h)})
and the  quantity  $C_2(r,\tau,\cL), \tau\in (0,1),\;\cL>0,$ is given in Section \ref{sec:constants}.
Its expression is rather cumbersome and it is why we do not present it right now.

\begin{theorem}
\label{th:deviation-l_p-gauss}
Let  $q\geq 1$, $p\geq 1$, $\tau\in (0,1)$, $\cL>0$  and $\cA\geq\mh^{-\frac{d}{2}}$ be fixed and let $\mathrm{H}$ be an arbitrary
countable subset of $\bH_d\big(\tau,\cL,\cA\big)$.

Then for any
 $\cA$, $\mh$ and $\tau$ satisfying (\ref{eq0:relations-between-parameters}) and
$K$ satisfying Assumption \ref{ass:kernel},
\begin{eqnarray*}
&&\bE\bigg\{\sup_{\vec{h}\in\mathrm{H}}
\Big[\big\|\xi_{\vec{h}}\big\|_p-\psi(\vec{h})\Big]_+\bigg\}^{q}
\leq  \bigg[C_4\cA e^{-e^{2\sqrt{2d|\ln(\mh)|}}}\bigg]^q, \quad\forall\mh\in \big(0,e^{-2}\big),
\end{eqnarray*}
where
$C_4$ depends on $\cK, p, q, b$ and $d$ only and its explicit expression can be found in Section \ref{sec:constants}.
\end{theorem}

Some  remarks are in order.

$\mathbf{1}^0.\;$
The statement of the theorem remains valid for any subset providing the measurability of  the corresponding  supremum. It explains, in particular,  why the upper function $\psi\big(\vec{h}\big)$ is independent of the choice of $\mathrm{H}$ and completely determined by the parameters
$\tau$, $\cL$ and $\cA$. It is worth noting that unlike Theorem \ref{th:deviation-l_p-gauss-general}  those proof is relatively standard the proof of Theorem  \ref{th:deviation-l_p-gauss} is rather long and  tricky.

\smallskip

$\mathbf{2}^0.\;$ Let us come back to the example of $\mathrm{H}$ (discussed after Theorem \ref{th:deviation-l_p-gauss-general}) consisting of "isotropic" constant functions. Obviously  $\mathrm{H}\subset\bH(\tau,\cL)$ for any $\tau\in (0,1)$ and $\cL\geq (2b)^{d\tau}$. Suppose additionally that there exist
$S\in\bN^*$ such that $h\in \big\{\mh e^{-s},\;s=0,\ldots, S \big\}$. Then, $\mathrm{H}\in\bB(\cA)$ whatever $\cA\geq \sqrt{\mh^{-d}e^{dS}}$
and $\bN^*_p\big(\vec{h},\cA\big)= \bN^*_p$ for any $\vec{h}\in\mathrm{H}$.

We deduce from Theorem \ref{th:deviation-l_p-gauss} that  in this case
$$
\psi(\vec{h})=h^{-\frac{d}{2}}\inf_{r\in\bN^*_p}C_2(r,\tau,\cL)
$$
and, therefore, $ \psi(\vec{h}) $  is sharp and better than the upper function found in Theorem \ref{th:deviation-l_p-gauss-general}.

\smallskip

$\mathbf{3}^0.\;$ It is impossible to compare both upper functions when an arbitrary subset of $\bH(\tau,\cL,\cA)$ is considered.  However they can be easily combined in such a way that
 the obtained  upper function is
smaller that both of them. Indeed, set $\Psi_\e\big(\vec{h}\big)=\psi_\e\big(\vec{h}\big)\wedge \psi\big(\vec{h}\big)$. First, we remark that
$$
\bigg\{\sup_{\vec{h}\in\mathrm{H}}
\Big[\big\|\xi_{\vec{h}}\big\|_p-\Psi_\e(\vec{h})\Big]_+\bigg\}^{q}\leq\bigg\{\sup_{\vec{h}\in\mathrm{H}}
\Big[\big\|\xi_{\vec{h}}\big\|_p-\psi_\e(\vec{h})\Big]_+\bigg\}^{q}+\bigg\{\sup_{\vec{h}\in\mathrm{H}}
\Big[\big\|\xi_{\vec{h}}\big\|_p-\psi(\vec{h})\Big]_+\bigg\}^{q}.
$$
Next, choose for instance
$
\mh=\mh_\e:=e^{-\sqrt{|\ln(\e)|}}, \;
\cA=\cA_\e:=e^{\ln^2(\e)}.
$
This yields
$$
\lim_{\e\to 0}\e^{-a}\cA e^{-e^{2\sqrt{2d|\ln(\mh)|}}}=0,\;\;\forall a>0,
$$
and moreover, for any $\tau\in (0,1)$ there exist $\e_0(\tau)$ such that for all $\e\leq \e_0(\tau) $ the relation (\ref{eq0:relations-between-parameters})
is fulfilled. In view of these remarks  we come to the following corollary of Theorem \ref{th:deviation-l_p-gauss}.

\begin{corollary}
\label{cor1:th:deviation-l_p-gauss}
Let assumptions of Theorem \ref{th:deviation-l_p-gauss} hold and let $\mh=\mh_\e$ and $\cA=\cA_\e$
Then  for any $\tau\in (0,1)$ and any $q\geq1$ one can find $\e(\tau,q)$ such that for any $\e\leq\e(\tau,q)$
\begin{eqnarray*}
&&\bE\bigg\{\sup_{\vec{h}\in\mathrm{H}}
\Big[\big\|\xi_{\vec{h}}\big\|_p-\Psi_\e\big(\vec{h}\big)\Big]_+\bigg\}^{q}
\leq \big\{(C_3+C_4)\e\big\}^q.
\end{eqnarray*}

\end{corollary}

$\mathbf{4}^0.\;$ Let $q=p$ and $\mathrm{H}$ consists of a single vector-function $\vec{\mathbf{h}}$. Then, obviously
$$
\bE\Big\{
\big\|\xi_{\vec{\mathbf{h}}}\big\|_p\Big\}^{p}= c_p \Big\|V^{-\frac{1}{2}}_{\vec{\mathbf{h}}}\Big\|^{p}_p,
$$
and the natural  question is:  may $\Big\{C\Big\|V^{-\frac{1}{2}}_{\vec{h}}\Big\|_{p},\;\vec{h}\in\mathrm{H}\Big\}$ be the upper function
over rather massive subset of $\bH_d\big(\tau,\cL,\cA\big)$ for some absolute constant $C$ ?

First, we note that $\Psi_\e\big(\vec{h}\big)\gg C\Big\|V^{-\frac{1}{2}}_{\vec{h}}\Big\|^{p}$ in general. Indeed, for any $ \vec{h}$
$$
\psi_\e(\vec{h})\geq C_1 \sqrt{|\ln(\e)|}\Big\|V^{-\frac{1}{2}}_{\vec{h}}\Big\|_{p}.
$$
Moreover, analyzing the expression of  $C_2(r,\tau,\cL)$  we easily seen that $C_2(r,\tau,\cL)\to\infty$ when $r\to\infty$ and, therefore,
$$
\psi_\e(\vec{h})=C_2(r(\vec{h}),\tau,\cL)\Big\|V^{-\frac{1}{2}}_{\vec{h}}\Big\|_{\frac{pr(\vec{h})}{r(\vec{h})-p}},
$$
where $r(\vec{h})<\infty$ for any $\vec{h}$. It remains to note that $\frac{pr}{r-p}>p$ for any $r<\infty$.

Next, it is interesting to note that unexpectedly  the negative answer on the aforementioned question comes from the lower bound for minimax risks over anisotropic Nikolskii classes proved in  \cite{lepski-kerk-08}. We have no place here to discuss this issue in detail
and only mention that if the answer would be positive it would contradict to the assertion of Theorem 2 in \cite{lepski-kerk-08}. In this context we conjecture that the upper function found in Corollary \ref{cor1:th:deviation-l_p-gauss} is sharp if an arbitrary subset of $\bH(\tau,\cL,\cA)$ is considered.

\subsection{Isotropic case} In this section we will suppose that $\vec{h}(\cdot)=\big(h(\cdot),\ldots,h(\cdot)\big)$ and consider the case $p\in [1,2]$.
We will show that under these restrictions the result similar to those obtained in Theorem \ref{th:deviation-l_p-gauss} can be proved \textsf{without} \textsf{any} \textsf{condition}
imposed on the set of bandwidths.

Note that in the isotropic case $V_{\vec{h}}(x)=h^{d}(x)$ and introduce
$$
\psi^*\big(\vec{h}\big)=\inf_{r\in\bN^*, r> d} C_2^*(r) \Big\|h^{-\frac{d}{2}}\Big\|_{p+\frac{1}{r}},
$$
where the explicit expression of $C_2^*(r)$ is given  in Section \ref{sec:constants3}.

\begin{theorem}
\label{th:deviation-l_p-gauss-isotropic}
Let  $q\geq 1$, $p\in [1,2]$,  be fixed, let $\mathrm{H}$ be an arbitrary
countable subset of $\mS_d(\mh)$
and $\vec{h}(\cdot)=\big(h(\cdot),\ldots,h(\cdot)\big)$ for any $\vec{h}\in\mathrm{H}$.
Suppose also that Assumption \ref{ass:kernel-new} is fulfilled.
Then,
\begin{eqnarray*}
&&\bE\bigg\{\sup_{\vec{h}\in\mathrm{H}}
\Big[\big\|\xi_{\vec{h}}\big\|_p-\psi^*\big(\vec{h}\big)\Big]_+\bigg\}^{q}
\leq \Big(C_5e^{\mh^{-d}}\Big)^{q},\quad \forall \mh\leq e^{-2},
\end{eqnarray*}
where $C_5$ depends on  $K, p, q, b$ and $d$ only and its explicit expression can be found in Section \ref{sec:constants3}.

\end{theorem}

 Coming back to the example of $\mathrm{H}$  consisting of constant functions we conclude that Theorem \ref{th:deviation-l_p-gauss-isotropic} generalizes  the result obtained in \cite{lepski-a}, Theorem 1,
 as well as the result given by Theorem \ref{th:deviation-l_p-gauss} when $p\in [1,2]$. Indeed,  we do not require here the finiteness of the set in which the bandwidth takes its values.

Although  the proof of the theorem is based upon the same approach, which is applied for proving  Theorem \ref{th:deviation-l_p-gauss}, it requires to use quite different arguments. Both assumptions isotropy and $p\in [1,2]$ are crucial for deriving the statement of Theorem \ref{th:deviation-l_p-gauss-isotropic}.

Combining the results of Theorems \ref{th:deviation-l_p-gauss-general} and \ref{th:deviation-l_p-gauss-isotropic} we arrive to the following assertion.

\begin{corollary}
\label{cor1:th:deviation-l_p-gauss-isotropic}
Let  $q\geq 1$, $p\in [1,2]$,  be fixed, let $\mathrm{H}$ be an arbitrary
countable subset of $\mS_d(\mh)$
and $\vec{h}(\cdot)=\big(h(\cdot),\ldots,h(\cdot)\big)$ for any $\vec{h}\in\mathrm{H}$.
Suppose also that Assumption \ref{ass:kernel-new} is fulfilled and choose $\mh=\mh_\e$. Then,
\begin{eqnarray*}
&&\bE\bigg\{\sup_{\vec{h}\in\mathrm{H}}
\Big[\big\|\xi_{\vec{h}}\big\|_p-\psi_e\big(\vec{h}\big)\wedge\psi^*\big(\vec{h}\big)\Big]_+\bigg\}^{q}
\leq \Big([C_3+C_5]\e\Big)^{q},\quad\forall \e\in \big(0,e^{-2}\big].
\end{eqnarray*}

\end{corollary}

\subsection{Example of the functional class $\bH(\tau,\cL,\cA)$}
\label{sec:examples}

For any  $\vec{\beta}\in (0,\infty)^d$, $\vec{r}\in [1,\infty]^d$ and $\vec{L}\in (0,\infty)^d$ let $\bN_d(\vec{\beta},\vec{r},\vec{L})$ denote anisotropic
Nikolskii class of functions on $\bR^d$, see \cite{Nikolski}, Chapter 4, Section 3.

Let  $\ell$ be an arbitrary integer number,
and let $w:\bR\to \bR$ be a compactly supported function satisfying
 $w\in\bC^{1}(\bR)$. Put
\begin{equation*}
\label{eq:w-function}
 w_\ell(y)=\sum_{i=1}^\ell \binom{\ell}{i} (-1)^{i+1}\frac{1}{i}w\Big(\frac{y}{i}\Big),\qquad
 K(t)=\prod_{j=1}^d w_\ell(t_j),\;\;\;\;t=(t_1,\ldots,t_d).
\end{equation*}
Although it will not be important for our  considerations here we note nevertheless that $K$ satisfies Assumption \ref{ass:kernel} with
$\cK=w_\ell$.

\smallskip

Let $\e,\mh\in\big(0,e^{-2}\big]$ be fixed and
set
$
\frac{1}{\beta}=\sum_{i=1}^d\frac{1}{\beta_i},\; \frac{1}{\upsilon}=\sum_{i=1}^d\frac{1}{r_i\beta_i}.
$
For any $j=1,\ldots,d$ let $S_\e(j)\in\bN^*$ be defined from the relation
\begin{equation}
\label{eq:example-def-S(j)}
 e^{-1}\e^{\frac{2\beta}{(2\beta+1)\beta_j}}<\mh e^{-S_\e(j)}\leq \e^{\frac{2\beta}{(2\beta+1)\beta_j}}.
\end{equation}
Without loss of generality we will assume that $\e$ is sufficiently small in order to provide the existence of $S_\e(j)$ for any $j$.
Put also
$$
\mH^{(j)}_\e=\{\mh_s=\mh e^{-s}, s\in\bN, \; s\geq S_\e(j)\},\qquad\mH_\e=\mH^{(1)}_\e\times\cdots\times\mH^{(d)}_\e
$$
and introduce for any $x\in (-b,b)^d$ and any $f\in\bN_d(\vec{\beta},\vec{r},\vec{L})$
$$
\vec{h}_{f}(x)=\arg\inf_{\vec{h}\in\mH_\e}\bigg[\Big|\int K_{\vec{h}}(t-x)f(t)\rd t-f(x)\Big|+\e V^{-\frac{1}{2}}_{\vec{h}}\bigg],\quad V_{\vec{h}}=\prod_{i=1}^d h_i.
$$
Define finally $\mathrm{H}=\left\{\vec{h}_{f},\;f\in\bN_d(\vec{\beta},\vec{r},\vec{L})\right\}$.

\begin{proposition}
\label{prop:set-of-bandw}
Let   $\vec{\beta}\in (0,\ell]^d$, $\vec{r}\in [1,p]^d$ and $\vec{L}\in (0,\infty)^d$ be given.

\par

\noindent 1)  For any $\tau\in(0,1)$ there exists $\cL>0$ such that
$$
\left\{\vec{h}_{f},\;f\in\bN_d(\vec{\beta},\vec{r},\vec{L})\right\}\subset\bH(\tau,\cL).
$$
2) If additionally $\upsilon(2+1/\beta)>p$ then there exists $C>0$ such that
$$
\left\{\vec{h}_{f},\;f\in\bN_d(\vec{\beta},\vec{r},\vec{L})\right\}\subset\bB\big(C\e^{-\frac{1}{2\beta+1}}\big).
$$
\end{proposition}
The explicit expression for the constants $\cL$ and $C$ can be found in the proof of the proposition which is postponed to Appendix.

The condition $\upsilon(2+1/\beta)>p$ appeared in the second assertion of the proposition is known as \textsf{the dense zone} in adaptive minimax estimation over the collection of anisotropic classes of smooth functions on $\bR^d$, see \cite{GL13}.

\section{Proof of Theorems \ref{th:deviation-l_p-gauss-general}--\ref{th:deviation-l_p-gauss-isotropic}}
\label{sec:proofs}

The proofs of these theorems are  based on several auxiliary
results, which for the citation convenience are formulated in Lemmas \ref{lem:talagrand-lifshits} and \ref{lem:birman-solomyak} below.

Furthermore, for any totaly bounded metric space $(\mathfrak{T},\varrho)$ we  denote by $\mE_{\varrho,\mathfrak{T}}(\delta),\;\delta>0$, the $\delta$-entropy of $\mT$ measured in $\varrho$, i.e. the logarithm of the minimal number of  $\varrho$-balls of radius $\delta>0$ needed to cover $\mT$.

$\mathbf{1^0.}\;$ The results formulated in   Lemma \ref{lem:talagrand-lifshits} can be found in  \cite{Talagrand}, Proposition 2.2, and   \cite{Lif}, Theorems 14.1 and 15.2.
\begin{lemma}
\label{lem:talagrand-lifshits}
Let $(Z_t,\; t\in\bT)$ be a centered, bounded on $\bT$, gaussian random function.

 I) For any $u>0$
$$
\bP\bigg\{\sup_{t\in\bT}Z_t\geq \bE\Big(\sup_{t\in\bT}Z_t\Big)+u \bigg\}\leq e^{-\frac{u^{2}}{2\sigma^2}},
$$
where $\sigma^2=\sup_{t\in\bT}\bE\big( Z^2_t\big).$

 II) Let $\bT$ be equipped with intrinsic semi-metric $\rho(t,t^\prime):=\sqrt{\bE\left(Z_{t}-Z_{t^\prime}\right)},\;t,t^\prime\in\bT$.
 Then
$$
\bE\Big(\sup_{t\in\bT}Z_t\Big)\leq D_{\bT,\rho}:=4\sqrt{2}\int_{0}^{\sigma/2}\sqrt{\mathfrak{E}_{\rho,\bT}(\delta)}\rd \delta.
$$
 III) If $D_{\bT,\rho}<\infty$ then the $(Z_t,\;t\in\bT)$ is bounded and uniformly continuous almost surely.

\end{lemma}

$\mathbf{2^0.}\;$ The  result formulated in Lemma \ref{lem:birman-solomyak} below is a particular case of Theorem  5.2 in  \cite{birman-sol}.

\smallskip

Let $\gamma>0, \gamma\notin\bN^*,$ $m\geq 1$ and $R> 0$ be fixed numbers and let $\Delta_k\subset\bR^{k},\;k\geq 1,$ be a given  cube with the sides parallel to the axis.
Recall that  $|y|$ denotes the euclidian norm of $y\in\bR^k$ and
  $\lfloor\gamma\rfloor$ is the  integer part of $\gamma$.
Set also $D^{\mathbf{n}}=\frac{\partial^{|\mathbf{n}|}}{\partial y_1^{n_1}\cdots\partial y_k^{n_k}},\; \mathbf{n}=(n_1,\ldots,n_k)\in\bN^k$.

\par

Denote by $\bS^{\gamma}_m\big(\Delta_k\big)$  the Sobolev-Slobodetskii space, i.e. the set of functions  $F:\Delta_k\to\bR$ equipped  with the norm
$$
\|F\|_{\gamma,m}=\bigg(\int_{\Delta_k}\big|F(y)\big|^m\rd y\bigg)^{\frac{1}{m}}+
\bigg(\sum_{|\mathbf{n}|=\lfloor\gamma\rfloor}\int_{\Delta_k}\int_{\Delta_k}
\frac{\big|D^{\mathbf{n}}F(y)-D^{\mathbf{n}}F(z)\big|^m}{|y-z|^{k+m(\gamma-\lfloor\gamma\rfloor)}}
\rd y\rd z\bigg)^{\frac{1}{m}}.
$$
Denote by $\bS^{\gamma}_m\big(\Delta_k,R\big)=\left\{F:\Delta_k\to\bR:\; \|F\|_{\gamma,m}\leq R\right\}$ the ball of radius $R$ in this space and set
$$
\lambda_k\big(\gamma,m,R,\Delta_k\big)=\inf\left\{c:\;\; \sup_{\delta\in (0,R]}\delta^{k/\gamma}\mE_{\|\cdot\|_{2},\bS^{\gamma}_m\big(\Delta_k,R\big)}(\delta)\leq c\right\}.
$$

\begin{lemma}
\label{lem:birman-solomyak} $\lambda_k\big(\gamma,m,1,\Delta_k\big)<\infty$
for any bounded $\Delta_k$ and $\gamma, m,k$ satisfying $\gamma>k/m-k/2$.

\end{lemma}
In view of the obvious relation
$
\mE_{\|\cdot\|_{2},\bS^{\gamma}_m\big(\Delta_k,R\big)}(\delta)=\mE_{\|\cdot\|_{2},\bS^{\gamma}_m\big(\Delta_k,1\big)}(\delta/R)
$
one has for any $R>0$
\begin{equation}
\label{eq:cor-birman-solomjak}
\lambda_k\big(\gamma,m,R,\Delta_k\big)= R^{k/\gamma}\lambda\big(\gamma,m,1,\Delta_k\big).
\end{equation}

\subsection{Proof of  Theorem \ref{th:deviation-l_p-gauss-general}}

For any multi-index $\mathbf{s}\in\bN^d$ set $\vec{\mh}_{\mathbf{s}}=\big(\mh_{s_1},\ldots,\mh_{s_d}\big)$, $V_{\mathbf{s}}=\prod_{j=1}^d \mh_{\mathbf{s}_j}$ and introduce
$$
\eta_{\mathbf{s}}(x)=\Big(V_{\mathbf{s}}\Big)^{-\frac{1}{2}}\int K_{\vec{\mh}_{\mathbf{s}}}\left(t-x\right)W(\rd t),
\qquad \eta_{\mathbf{s}}=\Big(\big|\ln{\big(\e V_{\mathbf{s}}\big)}\big|\Big)^{-\frac{1}{2}}\sup_{x\in (-b,b)^d}\big|\eta_{\mathbf{s}}(x)\big|.
$$
Note first that
for  any  $\vec{h}\in\mathrm{H}$ and any $\mathbf{s}\in\bN^d$   we obviously have
$$
\big|\xi_{\vec{h}}(x)\big|\leq \eta_{\mathbf{s}}  V^{-\frac{1}{2}}_{\mathbf{s}}\sqrt{\big|\ln{\big(\e V_{\mathbf{s}}\big)}\big|}, \quad \forall x\in\Lambda_{\mathbf{s}}\big[\vec{h}\big],
$$
and, therefore,
\begin{eqnarray*}
&& \Big\|\xi_{\vec{h}} \Big\|^p_{p}\leq \sum_{\mathbf{s}\in\bN^d}\eta^p_{\mathbf{s}}\Big(\big|\ln{\big(\e V_{\mathbf{s}}\big)}\big|V^{-1}_{\mathbf{s}}\Big)^{\frac{p}{2}}\nu_d\big(\Lambda_{\mathbf{s}}\big[\vec{h}\big]\big).
\end{eqnarray*}
Since
$$
\Big\| V^{-\frac{1}{2}}_{\vec{h}}\sqrt{\big|\ln{\big(\e V_{\vec{h}}\big)}\big|}\Big\|^p_{p}=\sum_{\mathbf{s}\in\bN^d}\Big(\big|\ln{\big(\e V_{\mathbf{s}}\big)}\big|V^{-1}_{\mathbf{s}}\Big)^{\frac{p}{2}}\nu_d\big(\Lambda_{\mathbf{s}}\big[\vec{h}\big]\big),
$$
using obvious inequality $\big(y^{1/p}-z^{1/p}\big)_+\leq \big[(y-z)_+\big]^{1/p},\; y,z\geq 0, p\geq 1,$ we obtain for any $\vec{h}\in\mathrm{H}$
$$
\bigg(\Big\|\xi_{\vec{h}} \Big\|_{p}-\psi_{\e}(\vec{h})\bigg)_+\leq
\bigg[\sum_{\mathbf{s}\in\bN^d}\Big(\big|\ln{\big(\e V_{\mathbf{s}}\big)}\big|V^{-1}_{\mathbf{s}}\Big)^{\frac{p}{2}}\Big(\eta^p_{\mathbf{s}}-C_1\Big)_+\bigg]^{\frac{1}{p}}.
$$
Noting that the right hand side of the latter inequality is independent of $\vec{h}$ and denoting $\tilde{q}=(q/p)\vee 1$
we obtain using Jensen and triangle inequalities
\begin{eqnarray}
\label{eq00001:proof-th:deviation-l_p-gauss}
&&
\bE\bigg\{\sup_{\vec{h}\in\mathrm{H}}
\Big[\big\|\xi_{\vec{h}}\big\|_p-\psi_{\e}(\vec{h})\Big]_+\bigg\}^{q}
\leq\bigg[\sum_{\mathbf{s}\in\bN^d}\Big(\big|\ln{\big(\e V_{\mathbf{s}}\big)}\big|V^{-1}_{\mathbf{s}}\Big)^{\frac{p}{2}}
\Big\{\bE\big(\eta^p_{\mathbf{s}}-C_1\big)^{\tilde{q}}_+\Big\}^{\frac{1}{\tilde{q}}}\bigg]^{\frac{q}{p}}.
\end{eqnarray}

Let $\mathbf{s}\in\bN^d$ be fixed. We have
\begin{eqnarray}
\label{eq000400:proof-th:deviation-l_p-gauss}
\bE\big(\eta^p_{\mathbf{s}}-C_1\big)^{\tilde{q}}_+&=&\tilde{q}\int_0^\infty z^{\tilde{q}-1}\bP\big\{\eta^p_{\mathbf{s}}\geq C_1 +z\big\}\rd z
\nonumber
\\
&=& \tilde{q}\int_0^\infty z^{\tilde{q}-1}\bP\bigg\{\sup_{x\in (-b,b)^d}\big|\eta_{\mathbf{s}}(x)\big|\geq \big[C_1 +z\big]^{\frac{1}{p}}\sqrt{\big|\ln{\big(\e V_{\mathbf{s}}\big)}\big|}\bigg\}\rd z.
\end{eqnarray}
Set $\mathfrak{z}=\big[C_1 +z\big]^{\frac{1}{p}}\sqrt{\big|\ln{\big(\e V_{\mathbf{s}}\big)}\big|}$
and  prove that
\begin{eqnarray}
\label{eq020:proof-th:deviation-l_p-gauss}
\bP\Big\{\sup_{x\in (-b,b)^d}\big|\eta_{\mathbf{s}}(x)\big|\geq \mathfrak{z}\Big\}\leq 2\big(\e V_{\mathbf{s}}\big)^{2(q\vee p)}
\exp\bigg(-\frac{z^{\frac{2}{p}}}{8\|\cK\|^{2d}_2}\bigg),\quad \forall z\geq 0.
\end{eqnarray}

\par

Since $\eta_{\mathbf{s}}(\cdot)$ is zero mean gaussian random field we get in
 view of obvious relation $\sup_x|\eta_{\mathbf{s}}(x)|=[\sup_x \eta_{\mathbf{s}}(x)]\vee [\sup_x \{-\eta_{\mathbf{s}}(x)\}]$
\begin{eqnarray}
\label{eq4:proof-th:deviation-l_p-gauss}
&&\bP\bigg\{\sup_{x\in(-b,b)^d}\big|\eta_{\mathbf{s}}(x)\big|\geq \mathfrak{z}\bigg\}\leq 2\bP\bigg\{\sup_{x\in (-b,b)^d}\eta_{\mathbf{s}}(x)\geq \mathfrak{z}\bigg\}.
\end{eqnarray}
Let $\rho$ denote intrinsic semi-metric of $\eta_{\mathbf{s}}(\cdot)$ on $(-b,b)^d$.

We have for any $x,x^\prime\in(-b,b)^d$ in view of Assumption \ref{ass:kernel-new-new}
\begin{eqnarray}
\label{eq400:proof-th:deviation-l_p-gauss}
\rho^{2}\big(x,x^\prime\big)&\leq&\int\Big[K(u)-K\big(\vec{\mh}_{\mathbf{s}}^{-1}(x-x^\prime)+u\big)\Big]^2\rd u
\nonumber\\
&=&2\|K\|^2_2-2\int_{[-a,a]^d}K(u)K\big(\vec{\mh}_{\mathbf{s}}^{-1}(x-x^\prime)+u\big)\rd u
\nonumber\\
&=&-2\int_{[-a,a]^d}K(u)\Big[K\big(\vec{\mh}_{\mathbf{s}}^{-1}(x-x^\prime)+u\big)-K(u)\Big]\rd u
\nonumber\\
&\leq&
2L\|K\|_2\big|\vec{\mh}_{\mathbf{s}}^{-1}(x-x^\prime)\big|\leq 2L\|K\|_2V_{\mathbf{s}}^{-1}\big|x-x^\prime\big|.
\end{eqnarray}
It yields  for any $\delta>0$
\begin{eqnarray}
\label{eq401:proof-th:deviation-l_p-gauss}
\mE_{\rho,(-b,b)^d}(\delta)\leq dc_1+d\big|\ln\big(V_{\mathbf{s}}\big)\big|+2d\big[\ln(1/\delta)\big]_{+},
\end{eqnarray}
where $c_1=\big|\ln{\big(4bL\|K\|_2\big)}\big|$.

Note that
$
\sigma^{2}:=\sup_{x\in(-b,b)^d}\bE\left(\eta^2_{\mathbf{s}}(x)\right)= \|K\|^2_2
$
and, therefore,
\begin{eqnarray}
\label{eq5:proof-th:deviation-l_p-gauss}
D_{(-b,b)^d,\rho}\leq \sqrt{d}\left(c_2+2\sqrt{2}\|K\|_2\sqrt{\big|\ln\big(V_{\mathbf{s}}\big)\big|}\right),
\end{eqnarray}
where
$
c_2=2\|K\|_2\sqrt{2c_1}+4\sqrt{2}\int_{0}^{2^{-1}\|K\|_2}\sqrt{\big[\ln(1/\delta)]_+}\rd \delta.
$

Thus, using the second assertion of Lemma \ref{lem:talagrand-lifshits} we have
$$
\mathbf{E}:=\bE\Big(\sup_{x\in(-b,b)^d}\eta_{\mathbf{s}}(x)\Big)\leq 2\sqrt{2d\pi}+ 2\sqrt{2dc_1}\|K\|_2 +2\sqrt{2d}\|K\|_2\sqrt{\big|\ln\big(V_{\mathbf{s}}\big)\big|}.
$$
Here we have  used that
$
4\sqrt{2}\int_{0}^{2^{-1}\|K\|_2}\sqrt{\big[\ln(1/\delta)]_+}\rd \delta\leq 2\sqrt{2\pi}.
$

Note that in view of the definition of $C_1$
$$
\mathfrak{z}-\mathbf{E}\geq 2^{-1}C^{\frac{1}{p}}_1\sqrt{\big|\ln\big(\e V_{\mathbf{s}}\big)\big|}-\mathbf{E}+2^{-1}z^{\frac{1}{p}}\geq
2\sqrt{(q\vee p)}\|K\|_2\sqrt{\big|\ln\big(\e V_{\mathbf{s}}\big)\big|}+2^{-1}z^{\frac{1}{p}}.
$$
Remark that  the third  assertion  of Lemma \ref{lem:talagrand-lifshits} and (\ref{eq5:proof-th:deviation-l_p-gauss})
implies that
the first assertion  of Lemma \ref{lem:talagrand-lifshits} is applicable with $\bT=(-b,b)^d$ and $Z_t=\eta_{\mathbf{s}}(x)$ and we get for any $\mathbf{s}\in\bN^d$
\begin{eqnarray*}
\bP\bigg\{\sup_{x\in(-b,b)^d}\eta_{\mathbf{s}}(x)\geq \mathfrak{z}\bigg\}\leq
\big(\e V_{\mathbf{s}}\big)^{2(q\vee p)}
\exp\bigg(-\frac{z^{\frac{2}{p}}}{8\|K\|^{2}_2}\bigg).
\end{eqnarray*}
Thus, the inequality
 (\ref{eq020:proof-th:deviation-l_p-gauss}) follows now from (\ref{eq4:proof-th:deviation-l_p-gauss}).
 We obtain from  (\ref{eq000400:proof-th:deviation-l_p-gauss}) and (\ref{eq020:proof-th:deviation-l_p-gauss})
\begin{eqnarray}
\label{eq2:proof-th:deviation-l_p-gauss}
&&\bE\big(\eta^p_{\mathbf{s}}-C_1\big)^{\tilde{q}}_+\leq
2\tilde{q}\big(\e V_{\mathbf{s}}\big)^{2(q\vee p)}
\int_0^\infty z^{\tilde{q}-1}\exp\bigg(-\frac{z^{\frac{2}{p}}}{8\|K\|^{2}_2}\bigg)\rd z=:c_3\big(\e V_{\mathbf{s}}\big)^{2(q\vee p)}.
\end{eqnarray}
Taking into account that
$
\big|\ln{\big(\e V_{\mathbf{s}}\big)}\big|\leq \big|\ln{(\e)}\big|V^{-1}_{\mathbf{s}},
$
since $\e,\mh\leq e^{-2}$,
we deduce from  (\ref{eq00001:proof-th:deviation-l_p-gauss}) and (\ref{eq2:proof-th:deviation-l_p-gauss}) that
\begin{eqnarray*}
\label{eq07:proof-th:deviation-l_p-gauss}
&&\bE\bigg\{\sup_{\vec{h}\in\mathrm{H}}
\Big[\big\|\xi_{\vec{h}}\big\|_p-\psi_{\e}(\vec{h})\Big]_+\bigg\}^{q}\leq (c_3)^{\frac{q}{\tilde{q}p}}\e^q
\bigg[\sum_{\mathbf{s}\in\bN^d}V_{\mathbf{s}}^{p}\bigg]^{\frac{q}{p}}\leq 2^{\frac{dq}{p}}(c_3)^{\frac{q}{\tilde{q}p}}\e^q=(C_3\e)^q.
\end{eqnarray*}

\epr

\subsection{Proof of Theorem \ref{th:deviation-l_p-gauss}}

\subsubsection{Auxiliary lemma}

Set
$
\lambda^{*}(\gamma,m)=\lambda_1\big(\gamma,m,1,[-a-b,a+b]\big),
$
 where recall the number $a>0$ is involved in Assumption \ref{ass:kernel} and $\lambda_k(\cdot,\cdot,\cdot,\cdot),\;k\in\bN^d$ is defined in
 Lemma \ref{lem:birman-solomyak}.

If $d\geq 2$ denote $\mathrm{x}=(x_2,\ldots,x_d)$ and define for any $\vec{\eta}\in\mathrm{H}$ and any $\mathrm{x}\in(-b,b)^{d-1}$
$$
\lambda_{\vec{\eta},\mathbf{s}}(\mathrm{x})=\bigg[\int_{-b}^b\mathrm{1}_{\Lambda_{\mathbf{s}}[\vec{\eta}]}(x)\nu_1\big(\rd x_1\big)\bigg]^{\frac{\tau}{r}}.
$$
Later on for any $x\in(-b,b)^{d}$ we will use the following notation $x=(x_1,\mathrm{x})$. If $d=1$ the dependence of $\mathrm{x}$ should be omitted in all
formulas.
In particular,  if $d=1$ then
$
\lambda_{\eta_1,s_1}=\Big\{\nu_1\big(\Lambda_{s_1}[\eta_1]\big)\Big\}^{\frac{\tau}{r}}.
$

For any $\mathrm{x}\in(-b,b)^{d-1}$ and  $\mathbf{s}\in\bN^d$
introduce the set of functions $Q:\bR\to\bR$
$$
\cQ_{\mathrm{x},\mathbf{s}}=\left\{ Q(\cdot)=\lambda^{-1}_{\vec{\eta},\mathbf{s}}(\mathrm{x})\int_{-b}^b\mh_{s_1}^{-1/2}
\cK\bigg(\frac{\cdot-x_1}{\mh_{s_1}}\bigg)\ell(x_1)\mathrm{1}_{\Lambda_{\mathbf{s}}[\vec{\eta}]}(x_1,\mathrm{x})\nu_1(\rd x_1),\; \ell\in\bB_q, \;\vec{\eta}\in \mathrm{H}\right\}.
$$
where  $\bB_q=\left\{\ell:(-b,b)\to\bR:\;\; \int_{-b}^b |\ell(x_1)|^q\nu(\rd x_1)\leq 1\right\},\;1/q=1-1/r$.

If $\lambda_{\vec{\eta},\mathbf{s}}(\mathrm{x})=0$ put by continuity $Q\equiv 0$.
Put  finally  $\mu^{-1}=q^{-1}+\tau r^{-1}$ and note that $2>\mu>1$ since $\tau<1$ and $r>2$.
\begin{lemma}  For any $\mathrm{x}\in(-b,b)^{d-1}$, $\mathbf{s}\in\cS_d$ and any $\omega\in \big(1/\mu-1/2,1\big)$ one has
\label{lem:entropy-of-integral-operators}
\begin{eqnarray*}
\mE_{\|\cdot\|_2,\cQ_{\mathrm{x},\mathbf{s}}}(\epsilon)\leq
\lambda^*\big(\omega,\mu\big)R_\mu^{\frac{1}{\omega}} \;\mh_{s_1}^{\frac{1}{2\omega}-1}\epsilon^{-\frac{1}{\omega}}\;\quad \forall \epsilon\in
\Big(0,R_\mu\;\mh_{s_1}^{\frac{1}{2}-\omega}\Big],
\end{eqnarray*}
where
$
R_\mu=\left[\left\{2^{-1}\|\cK\|_{\frac{2\mu}{3\mu-2}}\right\}\vee\left\{\|\cK\|_1+2\left[5\big\{4L(a+1)\big\}^\mu+
4\big\{2\|\cK\|_1\big\}^\mu(2-\mu)^{-1}\right]^{\frac{1}{\mu}}\right\}\right].
$

\end{lemma}

\subsubsection{Constants and expressions}
\label{sec:constants}
Put
\begin{eqnarray*}
C_4&=&\bigg(\boldsymbol{\gamma}_{q+1}\sqrt{(\pi/2)}\big[1\vee(2b)^{qd}\big]\sum_{r\in\bN_p^*}e^{-e^{r}}\big[\big(r\sqrt{e}\big)^d
\|\cK\|_{\frac{2r}{r+2}}^d\big]^{\frac{q}{2}}\bigg)^{\frac{1}{q}},
\end{eqnarray*}
where $\boldsymbol{\gamma}_{q+1}$ is the $(q+1)$-th absolute moment of the standard normal distribution.

Introduce
$
\Omega=\Big\{\{\omega_1,\omega_2\}: \;\; \omega_1<1/2<\omega_2,\;\; [\omega_1,\omega_2]\subset (1/\mu-1/2,1)\Big\}
$
and set
\begin{eqnarray*}
&&C_2(r,\tau,\cL)=[1\vee (2b)^{d-1}\big]
\big[\cL^{\frac{1}{r}}+\cL^{\frac{\tau}{r}}\big(1-e^{-\frac{\tau p}{4}}\big)^{\frac{\tau-1}{r}}\big]
\big[\widetilde{C}_\mu+\widehat{C}\big]
+e^r\sqrt{2(1+q)}\big(r\sqrt{e}\big)^d\|\cK\|_{\frac{2r}{r+2}}^d;
\\*[2mm]
&&\widehat{C}_\mu=\bigg[\frac{r}{1-\tau}\int_{0}^\infty \big(u+\widetilde{C}_\mu\big)^{\frac{r+\tau-1}{1-\tau}}
 \exp{\Big\{-u^{2}\Big[2|\cK\|_2^{d-1}\|\cK\|_{\frac{2\mu}{3\mu-2}}\Big]^{-1}\Big\}}\rd u\bigg]^{\frac{1-\tau}{r}};
\\*[3mm]
&&\widetilde{C}_\mu=C_\mu +4^d\big(\sqrt{2e^{r}}+\sqrt{8\pi}\big)\|\cK\|_2^{d-1}\|\cK\|_{\frac{2\mu}{3\mu-2}};
\\*[2mm]
&&C_\mu=4\sqrt{2}\|\cK\|_2^{d-1}\inf_{\{\omega_1,\omega_2\}\in\Omega}\bigg[\sqrt{\lambda^*\big(\omega_2,\mu\big)}\big(1-[2\omega_2]^{-1}\big)R_\mu^{\frac{1}{2\omega_2}} +
\sqrt{\lambda^*\big(\omega_1,\mu\big)}\big([2\omega_1]^{-1}-1\big)R_\mu^{\frac{1}{2\omega_1}}\bigg].
\end{eqnarray*}

\subsubsection{Proof of Theorem \ref{th:deviation-l_p-gauss}} Put for brevity $C_2(r)=C_2(r,\tau,\cL)$ and let
$$
\psi_{r}(\vec{h})=C_2(r)\Big\| V^{-\frac{1}{2}}_{\vec{h}}\Big\|_{\frac{rp}{r-p}},\;r\in\bN^*_p.
$$
For any $\vec{h}\in\mathrm{H}$ define $r^*\big(\vec{h}\big)=\arg\inf_{r\in\bN^*_p(\vec{h},\cA)}\psi_r(\vec{h})$.
 Note that  $C_2(r)<\infty$ for any $r\in\bN^*_p$ and
$$
\psi_r(\vec{h})\geq C_2(r)\mh^{-d}\to \infty,\;\;r\to\infty,
$$
and, therefore, $r^*\big(\vec{h}\big)<\infty$ for any $\vec{h}\in\bB(\cA)$.
 The latter fact allows us to assert that
\begin{eqnarray}
\label{eq0002000000:proof-th:deviation-l_p-gauss}
\inf_{r\in\bN^*_p(\vec{h},\cA)}\psi_r(\vec{h})=\psi_{r^*\big(\vec{h}\big)}(\vec{h})=
:C_2\Big(r^*\big(\vec{h}\big),\tau,\cL\Big)\Big\|V^{-\frac{1}{2}}_{\vec{h}}\Big\|_{\frac{pr^*(\vec{h})}{r^*(\vec{h})-p}},
\end{eqnarray}
since $\bN_p^*(\vec{h},\cA)$ is a discrete set.

By definition $r^*\big(\vec{h}\big)\geq r_{\cA}\big(\vec{h}\big)$, where recall $r_{\cA}\big(\vec{h}\big)$ is defined in (\ref{eq:def-N_p(h)}).
Hence we get from  H\"older
inequality and  the definition of $r_{\cA}\big(\vec{h}\big)$
\begin{eqnarray}
\label{eq000200000:proof-th:deviation-l_p-gauss}
\Big\|V^{-\frac{1}{2}}_{\vec{h}}\Big\|_{\frac{pr^*(\vec{h})}{r^*(\vec{h})-p}}
\leq \big[1\vee(2b)^d\big]\Big\|V^{-\frac{1}{2}}_{\vec{h}}\Big\|_{\frac{pr_{\cA}(\vec{h})}{r_{\cA}(\vec{h})-p}}
\leq \cA \big[1\vee(2b)^d\big].
\end{eqnarray}
Set  for any $r\in\bN^*_p$ and $\vec{h}\in\mathrm{H}$
$$
\zeta_{\vec{h}}(r)=\Big\|V^{\frac{1}{2}}_{\vec{h}}\xi_{\vec{h}}\Big\|_{r},\qquad\zeta(r)=
\sup_{\vec{h}\in\mathrm{H}}\zeta_{\vec{h}}(r).
$$
We obtain for any $\vec{h}\in\mathrm{H}$, applying H\"older inequality
\begin{eqnarray}
\label{eq0002:proof-th:deviation-l_p-gauss}
\Big\|\xi_{\vec{h}}\Big\|_{p}&\leq& \inf_{r\in\bN^*_p}\bigg\{\zeta(r)\Big\|V^{-\frac{1}{2}}_{\vec{h}}\Big\|_{\frac{pr}{r-p}}\bigg\}
\leq \zeta\Big(r^*\big(\vec{h}\big)\Big)\Big\|V^{-\frac{1}{2}}_{\vec{h}}\Big\|_{\frac{pr^*(\vec{h})}{r^*(\vec{h})-p}}.
\end{eqnarray}
We deduce from    (\ref{eq0002000000:proof-th:deviation-l_p-gauss}),  (\ref{eq000200000:proof-th:deviation-l_p-gauss}) and   (\ref{eq0002:proof-th:deviation-l_p-gauss}) that for any $\vec{h}\in\mathrm{H}$
\begin{eqnarray*}
\label{eq0002001:proof-th:deviation-l_p-gauss}
&&\Big[\big\|\xi_{\vec{h}}\big\|_p-\inf_{r\in\bN_p(\vec{h})}\psi_r(\vec{h})\Big]^q_+
\leq \Big\|V^{-\frac{1}{2}}_{\vec{h}}\Big\|^q_{\frac{pr^*(\vec{h})}{r^*(\vec{h})-p}}
\bigg[\zeta\Big(r^*\big(\vec{h}\big)\Big)-C_2\Big(r^*\big(\vec{h}\big)\Big)\bigg]^q_+
\nonumber\\
&&\leq \cA^q\big[1\vee(2b)^{qd}\big]
\bigg[\zeta\Big(r^*\big(\vec{h}\big)\Big)-C_2\Big(r^*\big(\vec{h}\big)\Big)\bigg]^q_+
\leq \cA^q\big[1\vee(2b)^{qd}\big]
\sum_{r\in\bN^*_p}\big[\zeta(r)-C_2(r)\big]^q_+.
\end{eqnarray*}
To get the last inequality we have used that $r^*\big(\vec{h}\big)\in \bN^*_p$ for any $\vec{h}\in\mathrm{H}$.

Taking into account that the right hand side of the latter inequality is independent of $\vec{h}$
we get
\begin{eqnarray}
\label{eq00020022:proof-th:deviation-l_p-gauss}
&&\bE\bigg(\sup_{\vec{h}\in\mathrm{H}}
\Big[\big\|\xi_{\vec{h}}\big\|_p-\psi_r(\vec{h})\Big]_+\bigg)^q
\leq \big[1\vee(2b)^{qd}\big] \cA^q\sum_{r\in\bN^*_p}\bE\big[\zeta(r)-C_2(r)\big]^q_+.
\end{eqnarray}
We have also  for any $r\in\bN^*_p$
\begin{eqnarray}
\label{eq0002002:proof-th:deviation-l_p-gauss}
\bE\big[\zeta(r)-C_2(r)\big]^q_+=q\int_{0}^{\infty}z^{q-1}\bP\Big\{\zeta(r)\geq C_2(r)+z\Big\}\rd z.
\end{eqnarray}


$1^0.\;$ Our goal now is to prove the following inequality: for any $z\geq 0$ and $r\in\bN^*_p$
\begin{eqnarray}
\label{eq6:proof-th:deviation-l_p-gauss}
\bP\Big\{\zeta(r)\geq C_2(r)+z\Big\}\leq
e^{-e^{r}}e^{-qe^{2\sqrt{2d|\ln(\mh)|}}}\;\exp{\Big\{-\Big(2\big(r\sqrt{e}\big)^d\|\cK\|_{\frac{2r}{r+2}}^d\Big)^{-1}z^2\Big\}}.
\end{eqnarray}
To do that we note first that in view of duality arguments
$$
\zeta(r)=\sup_{\vec{h}\in\mathrm{H}}\zeta_{\vec{h}}(r)=
\sup_{\vec{h}\in\mathrm{H}}\sup_{\vartheta\in\bB_{q,d}}\int_{(-b,b)^d} V^{\frac{1}{2}}_{\vec{h}}(x)\xi_{\vec{h}}(x)\vartheta(x)\nu_d(\rd x),
$$
where $\bB_{q,d}=\{\vartheta:(-b,b)^d\to\bR: \;\; \|\vartheta\|_q\leq 1\}$ and $1/q=1-1/r$.

Obviously
$$
\Upsilon_{\vec{h},\vartheta}:=\int_{(-b,b)^d} V^{\frac{1}{2}}_{\vec{h}}(x)\xi_{\vec{h}}(x)\vartheta(x)\nu_d(\rd x)
$$
is centered gaussian random function on $\mathrm{H}\times\bB_{q,d}$. Hence, if we show that for some $0<T<\infty$
\begin{eqnarray}
\label{eq600:proof-th:deviation-l_p-gauss}
\bE\big\{\zeta(r)\big\}\leq T,
\end{eqnarray}
then the first  assertion of Lemma \ref{lem:talagrand-lifshits} with
\begin{eqnarray}
\label{eq60001:proof-th:deviation-l_p-gauss}
\sigma_{\Upsilon}^{2}:=\sup_{\vec{h}\in\mathrm{H}}\sup_{\theta\in\bB_{q,d}}\bE\big\{\Upsilon_{\vec{h},\vartheta}\big\}^{2}
\end{eqnarray}
will be applicable to the random variable $\zeta(r)$.

\smallskip

$1^0\mathbf{a}.\;$ Let us bound from above $\sigma_{\Upsilon}$. By definition
$$
\Upsilon_{\vec{h},\vartheta}=\int\bigg[\int_{(-b,b)^d} V^{-\frac{1}{2}}_{\vec{h}}(x)K\bigg(\frac{t-x}{\vec{h}(x)}\bigg)\vartheta(x)\nu_d(\rd x)\bigg]W(\rd t)
$$
and, therefore,
$$
\sigma_{\Upsilon}=\sup_{\vec{h}\in\mathrm{H}}\sup_{\vartheta\in\bB_{q,d}}
\bigg[\int\bigg[\int_{(-b,b)^d} V^{-\frac{1}{2}}_{\vec{h}}(x)
K\bigg(\frac{t-x}{\vec{h}(x)}\bigg)\vartheta(x)\nu_d(\rd x)\bigg]^2\nu_d(\rd t)\bigg]^{\frac{1}{2}}.
$$
In view of triangle inequality and Assumption \ref{ass:kernel} ($\mathbf{ii}$)
$$
\sigma_{\Upsilon}\leq \sum_{\mathbf{s}\in\bN^d}\prod_{j=1}^{d}\mh^{-\frac{1}{2}}_{s_j}\sup_{\vartheta\in\bB_{q,d}}\bigg(\int\bigg[\int_{(-b,b)^d} \bigg|\prod_{j=1}^d\cK\bigg(\frac{t_j-x_j}{\mh_{s_j}}\bigg)\bigg|\big|\vartheta(x)\big|\nu_d(\rd x)\bigg]^2\nu_d(\rd t)\bigg)^{\frac{1}{2}}.
$$
Applying Young inequality and taking into account that $\vartheta\in\bB_{q,d}$ we obtain
\begin{eqnarray}
\label{eq145:proof-th:deviation-l_p-gauss}
\sigma_{\Upsilon}\leq \|\cK\|^{d}_{\frac{2r}{r+2}}\sum_{\mathbf{s}\in\bN^d}\prod_{j=1}^{d}\mh^{\frac{1}{r}}_{s_j}\leq
\big[1-e^{-\frac{1}{r}}\big]^{-d}\|\cK\|_{\frac{2r}{r+2}}^d\mh^{\frac{d}{r}}\leq \big(r\sqrt{e}\big)^d \|\cK\|_{\frac{2r}{r+2}}^d\mh^{\frac{d}{r}}.
\end{eqnarray}

$1^0\mathbf{b}.\;$ Let us prove (\ref{eq600:proof-th:deviation-l_p-gauss}).
Set for any $\mathbf{s}\in \bN^{d}$,
  and $\vec{h}\in\mathrm{H}$
$$
\xi_{\vec{h},\mathbf{s}}(x)=\mathrm{1}_{\Lambda_{\mathbf{s}}\big[\vec{h}\big]}(x)
\int \bigg[\prod_{i=1} ^d\mh^{-\frac{1}{2}}_{s_i}\cK\big((t_i-x_i)/\mh_{s_i}\big)\bigg]W(\rd t),\;\;
x\in (-b,b)^d.
$$
We obviously have for any $\vec{h}\in\mathrm{H}$
\begin{eqnarray}
\label{eq6001:proof-th:deviation-l_p-gauss}
\zeta^{r}_{\vec{h}}(r)=\Big\|V^{\frac{1}{2}}_{\vec{h}}\xi_{\vec{h}}\Big\|_{r}^r=\sum_{\mathbf{s}\in\bN^d}\big\|\xi_{\vec{h},\mathbf{s}}\big\|_{r}^r.
\end{eqnarray}
Moreover,  note that $\big|\xi_{\vec{h},\mathbf{s}}(x)\big|\leq \mathrm{1}_{\Lambda_{\mathbf{s}}\big[\vec{h}\big]}(x)\big|\ln\big(\e V_\mathbf{s}\big)\big|^{\frac{1}{2}}\eta_\mathbf{s}$ for any $x\in (-b,b)^{d}$, where, recall, $V_\mathbf{s}$ and $\eta_\mathbf{s}$ are defined in the beginning of  the proof of Theorem \ref{th:deviation-l_p-gauss-general}.
Since, we have proved that $\eta_\mathbf{s}$ is bounded almost surely, one gets
\begin{eqnarray}
\label{eq7002:proof-th:deviation-l_p-gauss}
\int_{-b}^b\big|\xi_{\vec{h},\mathbf{s}}(x)\big|^r\nu_1(\rd x_1)\leq \lambda^{r}_{\vec{h},\mathbf{s}}(\mathrm{x})\big|\ln\big(\e V_\mathbf{s}\big)\big|^{\frac{r}{2}}\eta_\mathbf{s}^{r}=0, \quad \text{if}\quad  \lambda_{\vec{h},\mathbf{s}}(\mathrm{x})=0.
\end{eqnarray}
On the other hand in view of duality arguments
\begin{eqnarray}
\label{eq601:proof-th:deviation-l_p-gauss}
\int_{-b}^b\big|\xi_{\vec{h},\mathbf{s}}(x)\big|^r\nu_1(\rd x_1)=\bigg[\sup_{\ell\in\bB_q}
\int_{-b}^b\xi_{\vec{h},\mathbf{s}}(x)\ell(x_1)\nu_1\big(\rd x_1\big)\bigg]^r,
\end{eqnarray}
where, recall,  $\bB_{q}=\left\{\ell:(-b,b)\to\bR:\;\; \int_{-b}^b |\ell(y)|^q\nu(\rd y)\leq 1\right\},\;1/q=1-1/r$.

\smallskip

Let $d\geq 2$.
Define for any  $\mathbf{s}\in\bN^d$ and  $\mathrm{x}\in(-b,b)^{d-1}$
\begin{eqnarray}
\label{eq8004:proof-th:deviation-l_p-gauss}
&&\varsigma_{\mathbf{s}}\big(Q,\mathrm{x}\big)=\int Q(t_1)G_{\mathbf{s}}(\mathrm{t},\mathrm{x})W(\rd t),\quad Q\in\cQ_{\mathrm{x},\mathbf{s}}.
\end{eqnarray}
 Here we have put $\mathrm{t}=(t_2,\ldots,t_d)$,  denoted $t=(t_1,\mathrm{t})$ for any $t\in\bR^d$, and set
$$
G_{\mathbf{s}}(\mathrm{t},\mathrm{x})=\prod_{i=2}^d \mh^{-\frac{1}{2}}_{s_i}\cK\big((\mathrm{t}_i-\mathrm{x}_i)/\mh_{s_i}\big),\;\; \mathrm{t}\in\bR^{d-1},\;\mathrm{x}\in(-b,b)^{d-1}.
$$
Below we will prove that $\varsigma_{\mathbf{s}}(\mathrm{x}):=\sup_{Q\in\cQ_{\mathbf{s},\mathrm{x}}}\varsigma_{\mathbf{s}}\big(Q,\mathrm{x}\big)$ is random variable. This is important  because its definition uses
 the supremum over $\cQ_{\mathbf{s},\mathrm{x}}$ which is not countable.

 \smallskip

 The following simple remark is crucial for all further consideration:  in view of (\ref{eq7002:proof-th:deviation-l_p-gauss}) and (\ref{eq601:proof-th:deviation-l_p-gauss}) for any  $\mathrm{x}\in(-b,b)^{d-1}$,  $\mathbf{s}\in\bN^d$   and  for any $\vec{h}\in\mathrm{H}$
\begin{eqnarray}
\label{eq7032:proof-th:deviation-l_p-gauss}
\int_{-b}^b\big|\xi_{\vec{h},\mathbf{s}}(x_1,\mathrm{x})\big|^r\nu_1(\rd x_1)\leq \lambda^r_{\vec{h},\mathbf{s}}(\mathrm{x})\varsigma_{\mathbf{s}}^r(\mathrm{x}).
\end{eqnarray}
Indeed,   if  $\lambda_{\vec{h},\mathbf{s}}(\mathrm{x})=0$  (\ref{eq7032:proof-th:deviation-l_p-gauss}) follows from (\ref{eq7002:proof-th:deviation-l_p-gauss}).
If $\lambda_{\vec{h},\mathbf{s}}(\mathrm{x})>0$ then
$$
\int_{-b}^b\xi_{\vec{h},\mathbf{s}}(x)\ell(x_1)\nu_1\big(\rd x_1\big)=\lambda_{\vec{h},\mathbf{s}}(\mathrm{x})\int Q(t_1)G_{\mathbf{s}}(\mathrm{t},\mathrm{x})W(\rd t),
$$
with
$
Q(\cdot)=\lambda^{-1}_{\vec{h},\mathbf{s}}(\mathrm{x})\int_{-b}^b\mh_{s_1}^{-1/2}
\cK\left(\frac{\cdot-x_1}{\mh_{s_1}}\right)\ell(x_1)\mathrm{1}_{\Lambda_{\mathbf{s}}[\vec{h}]}(x_1,\mathrm{x})\nu_1(\rd x_1)\in\cQ_{\mathrm{x},\mathbf{s}},
$
and (\ref{eq7032:proof-th:deviation-l_p-gauss}) follows from (\ref{eq601:proof-th:deviation-l_p-gauss}).
We get from (\ref{eq7032:proof-th:deviation-l_p-gauss}) for any  $\vec{h}\in\mathrm{H}$ and $\mathbf{s}\in\bN^*$  in view of Fubini theorem
\begin{eqnarray*}
\big\|\xi_{\vec{h},\mathbf{s}}\big\|_{r}^r&=&\int_{(-b,b)^{d-1}}\int_{b}^b\big|\xi_{\vec{h},\mathbf{s}}(x_1,\mathrm{x})\big|^r\nu_1(\rd x_1)\nu_{d-1}(\rd \mathrm{x})\leq \int_{(-b,b)^d}\lambda^r_{\vec{h},\mathbf{s}}(\mathrm{x})\varsigma_\mathbf{s}^r(\mathrm{x})\nu_{d-1}(\rd \mathrm{x})
\\
&=&\int_{(-b,b)^d}\varsigma_\mathbf{s}^r(\mathrm{x})
\bigg[\int_{-b}^b\mathrm{1}_{\Lambda_{\mathbf{s}}\big[\vec{h}\big]}(x)\nu_1\big(\rd x_1\big)\bigg]^{\tau}\nu_{d-1}(\rd \mathrm{x}).
\end{eqnarray*}
Taking into account that $\tau<1$ and applying H\"older inequality to the outer integral we get
\begin{eqnarray}
\label{eq6010:proof-th:deviation-l_p-gauss}
&&\big\|\xi_{\vec{h},\mathbf{s}}\big\|_{r}^r
\leq\nu^\tau_d\Big(\Lambda_{\mathbf{s}}\big[\vec{h}\big]\Big)\bigg\{\int_{(-b,b)^d}
\varsigma_\mathbf{s}^{\frac{r}{1-\tau}}(\mathrm{x})\nu_{d-1}(\rd \mathrm{x})\bigg\}^{1-\tau},\quad\forall \mathbf{s}\in\bN^*.
\end{eqnarray}
If $d=1$ putting  $G_{\mathbf{s}}(\mathrm{t},\mathrm{x})\equiv 1$ in (\ref{eq8004:proof-th:deviation-l_p-gauss}),  we obtain using the same arguments
\begin{eqnarray}
\label{eq6020:proof-th:deviation-l_p-gauss}
\big\|\xi_{h_1,s_1}\big\|_{r}^r\leq \nu^\tau_d\Big(\Lambda_{s_1}\big[h_1\big]\Big)\varsigma_{s_1},\qquad
\varsigma_{s_1}=\sup_{Q\in\cQ_{s_1}}\varsigma_{s_1}\big(Q\big).
\end{eqnarray}

$1^0\mathbf{b1}.\;$ Let us prove some bounds used in the sequel. Let $S\in\bN$ be the number satisfying
$
e^{-1}<\mh^{d}e^{-S}\cA^4\leq 1,
$
and set $\cS_d=\{0,1,\ldots S\}^d$ and $\bar{\cS_d}=\bN^d\setminus\cS_d$. If such $S$ does not exist we will assume that $\cS_d=\emptyset$ and later on the supremum
over empty set is assumed to be $0$.

Set also $\cS^*_d=\{\mathbf{s}\in\bN^d:\; \cA^4V_\mathbf{s}\leq 1\}$.
Note that
$
V_\mathbf{s}\leq \mh^{d}e^{-S}\leq\cA^{-4}
$
for any $\mathbf{s}\in\bar{\cS_d}$ and, therefore,
\begin{eqnarray}
\label{eq6011:proof-th:deviation-l_p-gauss}
\bN^d\setminus\cS^*_d\subseteq \cS_d.
\end{eqnarray}
Putting for brevity $\mathrm{r}=r_{\cA}(\vec{h})$, we have for any $\mathbf{s}\in\bN^d$ and any $\vec{h}\in\bB(\cA)$
$$
\big(V_{\mathbf{s}}\big)^{-\frac{p\mathrm{r}}{2(\mathrm{r}-p)}}\nu_d\Big(\Lambda_{\mathbf{s}}\big[\vec{h}\big]\Big)\leq \sum_{\mathbf{k}\in\bN^d}
\big(V_{\mathbf{k}}\big)^{-\frac{p\mathrm{r}}{2(\mathrm{r}-p)}}\nu_d\Big(\Lambda_{\mathbf{k}}\big[\vec{h}\big]\Big)=
\Big\| V^{-\frac{1}{2}}_{\vec{h}}\Big\|_{\frac{p\mathrm{r}}{\mathrm{r}-p}}^{\frac{p\mathrm{r}}{\mathrm{r}-p}}
\leq \cA^{\frac{p\mathrm{r}}{\mathrm{r}-p}}.
$$
The last inequality follows from the definition of $r_{\cA}(\vec{h})$.

Taking into account that $\frac{p\mathrm{r}}{\mathrm{r}-p}>p$ and that $V_\mathbf{s}<1$ we get in view of the definition of $\cS^*_d$
\begin{eqnarray}
\label{eq6012:proof-th:deviation-l_p-gauss}
\nu_d\Big(\Lambda_{\mathbf{s}}\big[\vec{h}\big]\Big)\leq V_{\mathbf{s}}^{\frac{p}{4}},\qquad\forall \vec{h}\in\bB(\cA),\;\;\forall\mathbf{s}\in\cS^*_d.
\end{eqnarray}

$1^0\mathbf{b2}.\;$ Set $\varsigma(\mathrm{x)}=\sup_{\mathbf{s}\in\cS_d}\varsigma_\mathbf{s}(\mathrm{x})$ and let $d\geq 2$.

We have in view of (\ref{eq6001:proof-th:deviation-l_p-gauss}), (\ref{eq6010:proof-th:deviation-l_p-gauss}),
(\ref{eq6011:proof-th:deviation-l_p-gauss}) and (\ref{eq6012:proof-th:deviation-l_p-gauss}) for any $\vec{h}\in\mathrm{H}$
\begin{eqnarray*}
\zeta^r_{\vec{h}}(r)
&\leq& \cL\bigg\{\int_{(-b,b)^{d-1}}\varsigma^{\frac{r}{1-\tau}}(\mathrm{x})\nu_{d-1}(\rd \mathrm{x})\bigg\}^{1-\tau}
\\
&\qquad+&\sum_{\mathbf{s}\in\cS^*_d}\nu^{\tau^2}_d\Big(\Lambda_{\mathbf{s}}\big[\vec{h}\big]\Big)
 V_{\mathbf{s}}^{\frac{\tau(1-\tau) p}{4}}\bigg\{\int_{(-b,b)^{d-1}}
\varsigma_\mathbf{s}^{\frac{r}{1-\tau}}(\mathrm{x})\nu_{d-1}(\rd \mathrm{x})\bigg\}^{1-\tau}.
\end{eqnarray*}
Here we have also used that $\mathrm{H}\subset\bH_d(\tau,\cL)$.

Applying H\"older  inequality with exponents $1/\tau$ and $1/(1-\tau)$
to the sum appeared in the second term in the right hand side of the latter inequality we get
\begin{eqnarray*}
&&\sum_{\mathbf{s}\in\cS^*_d}\nu^{\tau^2}_d\Big(\Lambda_{\mathbf{s}}\big[\vec{h}\big]\Big)
 V_{\mathbf{s}}^{\frac{\tau(1-\tau) p}{4}}\bigg\{\int_{(-b,b)^{d-1}}
\varsigma_\mathbf{s}^{\frac{r}{1-\tau}}(\mathrm{x})\nu_{d-1}(\rd \mathrm{x})\bigg\}^{1-\tau}
\\
&&
\leq \cL^{\tau}\bigg[ \sum_{\mathbf{s}\in\bN^d} V_{\mathbf{s}}^{\frac{\tau p}{4}}\int_{(-b,b)^{d-1}}
\varsigma_\mathbf{s}^{\frac{r}{1-\tau}}(\mathrm{x})\nu_{d-1}(\rd \mathrm{x})\bigg]^{1-\tau}.
\end{eqnarray*}
It yields for any  $\vec{h}\in\mathrm{H}$
\begin{eqnarray*}
\zeta^r_{\vec{h}}(r)
&\leq& \cL\bigg\{\int_{(-b,b)^{d-1}}\varsigma^{\frac{r}{1-\tau}}(\mathrm{x})\nu_{d-1}(\rd \mathrm{x})\bigg\}^{1-\tau}
\\
&\qquad+&\cL^{\tau}\bigg[ \sum_{\mathbf{s}\in\bN^d} V_{\mathbf{s}}^{\frac{\tau p}{4}}\int_{(-b,b)^{d-1}}
\varsigma_\mathbf{s}^{\frac{r}{1-\tau}}(\mathrm{x})\nu_{d-1}(\rd \mathrm{x})\bigg]^{1-\tau}.
\end{eqnarray*}

Noting that the right hand side of the obtained inequality is independent of $\vec{h}$ we get
\begin{eqnarray*}
\zeta(r)
&\leq& \cL^{\frac{1}{r}}\bigg\{\int_{(-b,b)^{d-1}}\varsigma^{\frac{r}{1-\tau}}(\mathrm{x})\nu_{d-1}(\rd \mathrm{x})\bigg\}^{\frac{1-\tau}{r}}
\\
&\qquad+&\cL^{\frac{\tau}{r}}\bigg[ \sum_{\mathbf{s}\in\bN^d} V_{\mathbf{s}}^{\frac{\tau p}{4}}\int_{(-b,b)^{d-1}}
\varsigma_\mathbf{s}^{\frac{r}{1-\tau}}(\mathrm{x})\nu_{d-1}(\rd \mathrm{x})\bigg]^{\frac{1-\tau}{r}}.
\end{eqnarray*}

Hence, applying Jensen inequality and Fubini theorem one has for any $d\geq 2$
\begin{eqnarray}
\label{eq602:proof-th:deviation-l_p-gauss}
\bE\big\{\zeta(r)\big\}&\leq&
\cL^{\frac{1}{r}}\bigg\{\int_{(-b,b)^{d-1}}\bE\Big(\varsigma^{\frac{r}{1-\tau}}(\mathrm{x})\Big)\nu_{d-1}(\rd \mathrm{x})\bigg\}^{\frac{1-\tau}{r}}
\nonumber\\
&\qquad+&\cL^{\frac{\tau}{r}}\bigg[ \sum_{\mathbf{s}\in\bN^d} V_{\mathbf{s}}^{\frac{\tau p}{4}}\int_{(-b,b)^{d-1}}
\bE\Big(\varsigma_\mathbf{s}^{\frac{r}{1-\tau}}(\mathrm{x})\Big)\nu_{d-1}(\rd \mathrm{x})\bigg]^{\frac{1-\tau}{r}}
\nonumber\\
&\leq& \cL^{\frac{1}{r}}\big[1\vee (2b)^{d-1}\big]\sup_{\mathrm{x}\in (-b,b)^{d-1}}
\bigg\{\bE\Big(\varsigma^{\frac{r}{1-\tau}}(\mathrm{x})\Big)\bigg\}^{\frac{1-\tau}{r}}
\nonumber\\
&\qquad+&\cL^{\frac{\tau}{r}}\big[1\vee (2b)^{d-1}\big]\big(1-e^{-\frac{\tau p}{4}}\big)^{\frac{\tau-1}{r}}\sup_{\mathbf{s}\in\bN^d}
\sup_{\mathrm{x}\in (-b,b)^{d-1}}
\bigg\{\bE\Big(\varsigma_s^{\frac{r}{1-\tau}}(\mathrm{x})\Big)\bigg\}^{\frac{1-\tau}{r}}.
\end{eqnarray}
Here we have also used that $(1-\tau)/r<1$.

If $d=1$ repeating  previous computations we obtain from  (\ref{eq6001:proof-th:deviation-l_p-gauss}) and (\ref{eq6020:proof-th:deviation-l_p-gauss})
\begin{eqnarray}
\label{eq60201:proof-th:deviation-l_p-gauss}
\bE\big\{\zeta(r)\big\}\leq
\cL^{\frac{1}{r}}\bE\varsigma
+\cL^{\frac{\tau}{r}}\big(1-e^{-\frac{\tau p}{4}}\big)^{\frac{\tau-1}{r}} \sup_{s\in\bN}\Big[
\bE\Big(\varsigma_s^{\frac{r}{1-\tau}}\Big)\Big]^{\frac{1-\tau}{r}}.
\end{eqnarray}
In  what follows $\mathrm{x}$ is assumed to be fixed that allows us not to separate cases $d=1$ and $d\geq 2$.
\smallskip

$1^0\mathbf{b3}.\;$ Let $\mathrm{x}\in(-b,b)^{d-1}$  be fixed. First let us bound from above
$$
\bE\varsigma_s(\mathrm{x}):=\bE\Big\{\sup_{Q\in\cQ_{\mathbf{s},\mathrm{x}}}\varsigma_{\mathbf{s}}\big(Q,\mathrm{x}\big)\Big\},\;\;
 \mathbf{s}\in\bN^d,\qquad
\bE\varsigma(\mathrm{x}):=\bE\Big\{\sup_{\mathbf{s}\in\cS_d}\sup_{Q\in\cQ_{\mathbf{s},\mathrm{x}}}\varsigma_{\mathbf{s}}\big(Q,\mathrm{x}\big)\Big\}
$$
Note that $\varsigma_{\mathbf{s}}\big(Q,\mathrm{x}\big)$ is zero-mean gaussian random function on $\cQ_{\mathbf{s},\mathrm{x}}$.
Our objective now is to show that the  assertion II of Lemma \ref{lem:talagrand-lifshits} is applicable with
$Z_t=\varsigma_{\mathbf{s}}\big(Q,\mathrm{x}\big), t=Q,$ and $\bT=\cQ_{\mathbf{s},\mathrm{x}}$.

Note that the intrinsic semi-metric of $\varsigma_{\mathbf{s}}\big(Q,\mathrm{x}\big)$ is given by
$$
\rho^{2}\big(Q,\widetilde{Q}\big)=\int G^{2}_{\mathbf{s}}\big(\mathrm{t},\mathrm{x}\big)\left[Q(t_1)-\widetilde{Q}(t_1)\right]^{2}\nu_d(\rd t),\quad\forall Q,\widetilde{Q}\in\cQ_{\mathbf{s},\mathrm{x}},
$$
and, therefore,
$
\rho\big(Q,\widetilde{Q}\big)=\|\cK\|_2^{d-1}\|Q-\widetilde{Q}\|_2,
$
for any $Q,\widetilde{Q}\in\cQ_{\mathbf{s},\mathrm{x}}$.

Below we show that $\big(\cQ_{\mathbf{s},\mathrm{x}},\|\cdot\|_2\big)$ is totally bounded metric space and, moreover, the corresponding
 Dudley's integral is finite.
 The latter fact allows us to assert that $\varsigma_{\mathbf{s}}(\cdot,\mathrm{x})$ is almost surely continuous on
 $\cQ_{\mathbf{s},\mathrm{x}}$ that implies the measurability of $\varsigma_{\mathbf{s}}(\mathrm{x})$ as well as $\varsigma(\mathrm{x})$.

We obviously have
\begin{eqnarray}
\label{eq14:proof-th:deviation-l_p-gauss}
&& \mE_{\rho,\cQ_{\mathbf{s},\mathrm{x}}}(\delta)\leq
\mE_{||\cdot||_2,\cQ_{\mathbf{s},\mathrm{x}}}\big(\|\cK\|_2^{1-d}\delta\big),\;\; \forall \delta>0,
\end{eqnarray}
and, therefore,
\begin{eqnarray}
\label{eq141:proof-th:deviation-l_p-gauss}
&&  D_{\cQ_{\mathbf{s},\mathrm{x}},\rho}:=4\sqrt{2}\int_{0}^{2^{-1}\sigma_\mathbf{s}}\sqrt{\mathfrak{E}_{\rho,\cQ_{\mathbf{s},\mathrm{x}}}(\delta)}\rd \delta\leq4\sqrt{2}\|\cK\|_2^{d-1}\int_{0}^{\tilde{\sigma}_\mathbf{s}}\sqrt{\mathfrak{E}_{||\cdot||_2,\cQ_{\mathbf{s},\mathrm{x}}}(\delta)}\rd \delta,
\end{eqnarray}
where $\tilde{\sigma}_\mathbf{s}=2^{-1}\sigma_\mathbf{s}\|\cK\|_2^{1-d}$ and
$$
\sigma_\mathbf{s}:=\Big[\sup_{Q\in\cQ_{\mathbf{s},\mathrm{x}}}
\bE\big\{\varsigma^2_{\mathbf{s}}\big(Q,\mathrm{x}\big)\big\}\Big]^{\frac{1}{2}}=\|\cK\|_2^{d-1}\sup_{Q\in\cQ_{\mathbf{s},\mathrm{x}}}\|Q\|_2.
$$
We start with bounding from above the quantity $\sigma_\mathbf{s}$.

Recall that $\mu^{-1}=q^{-1}+\tau r^{-1}$. Applying Young inequality we have
$$
\|Q\|_2\leq\lambda^{-1}_{\vec{h},\mathbf{s}}(\mathrm{x})\mh_{s_1}^{1-\frac{1}{\mu}}
\bigg[\int_{-b}^b
\big|\ell(x_1)\big|^\mu\mathrm{1}_{\Lambda_{\mathbf{s}}\big[\vec{h}\big]}(x_1,\mathrm{x})\nu_1(\rd x_1)\bigg]^{\frac{1}{\mu}}\|\cK\|_{\frac{2\mu}{3\mu-2}}.
$$
Applying H\"older inequality to the integral in right hand side of the latter inequality and taking into account that
$\ell\in\bB_q$  we get
\begin{eqnarray}
\label{eq604:proof-th:deviation-l_p-gauss}
&&\bigg[\int_{-b}^b
\big|\ell(x_1)\big|^\mu\mathrm{1}_{\Lambda_{\mathbf{s}}\big[\vec{h}\big]}(x_1,\mathrm{x})\nu_1(\rd x_1)\bigg]^{\frac{1}{\mu}}\leq
\bigg[\int_{-b}^b
\mathrm{1}_{\Lambda_{\mathbf{s}}\big[\vec{h}\big]}(x_1,\mathrm{x})\nu_1(\rd x_1)\bigg]^{\frac{1}{\mu}-\frac{1}{q}}
=\lambda_{\vec{h},\mathbf{s}}(\mathrm{x}).
\end{eqnarray}
Thus we obtain
\begin{eqnarray}
\label{eq605:proof-th:deviation-l_p-gauss}
&&\sigma_\mathbf{s}\leq \|\cK\|_2^{d-1}\|\cK\|_{\frac{2\mu}{3\mu-2}}\mh_{s_1}^{\frac{1-\tau}{r}}.
\end{eqnarray}
Putting $\sigma^*_\mathbf{s}=2^{-1}\|\cK\|_{\frac{2\mu}{3\mu-2}}\mh_{s_1}^{\frac{1-\tau}{r}}$ we deduce from (\ref{eq141:proof-th:deviation-l_p-gauss}) and (\ref{eq605:proof-th:deviation-l_p-gauss})
\begin{eqnarray}
\label{eq142:proof-th:deviation-l_p-gauss}
&&  D_{\cQ_{\mathbf{s},\mathrm{x}},\rho}\leq 4\sqrt{2}\|\cK\|_2^{d-1}\int_{0}^{\sigma^*_\mathbf{s}}\sqrt{\mathfrak{E}_{||\cdot||_2,\cQ_{\mathbf{s},\mathrm{x}}}(\delta)}\rd \delta.
\end{eqnarray}
Now let us bound from above $\bE\big\{\sup_{Q\in\cQ_{\mathbf{s},\mathrm{x}}}\varsigma_{\mathbf{s}}\big(Q,\mathrm{x}\big)\big\}$.

Recall that
$
\Omega=\Big\{\{\omega_1,\omega_2\}: \;\; \omega_1<1/2<\omega_2,\;\; [\omega_1,\omega_2]\subset (1/\mu-1/2,1)\Big\}.
$
Note that the condition $\omega_1>1/\mu-1/2$ implies
$
1/2-\omega_1<(1-\tau)r^{-1}
$
and, therefore
$$
\mh_{s_1}^{\frac{1-\tau}{r}}<\mh_{s_1}^{\frac{1}{2}-\omega_1}<\mh_{s_1}^{\frac{1}{2}-\omega_2}
$$
since $\mh_{s_1}\leq \mh\leq 1.$ It yields that $\big(0,\sigma^*_\mathbf{s}\big]\subset
\Big(0,R_\mu\;\mh_{s_1}^{\frac{1}{2}-\omega_1}\Big]\subset\Big(0,R_\mu\;\mh_{s_1}^{\frac{1}{2}-\omega_2}\Big]$, since $R_\mu\geq 2^{-1}\|\cK\|_{\frac{2\mu}{3\mu-2}}$.
Hence Lemma \ref{lem:entropy-of-integral-operators} is applicable to the computation of the integral in the right hand side of
(\ref{eq142:proof-th:deviation-l_p-gauss}).

Introduce the following notations:
$
A^{2}(\omega)=\lambda^*\big(\omega,\mu\big)R_\mu^{\frac{1}{\omega}} \;\mh_{s_1}^{\frac{1}{2\omega}-1}, \; \delta_0=\mh_{s_1}^{\frac{1}{2}}
$
and note that $\delta_0< \sigma^*_\mathbf{s}$.
 We get in view of Lemma \ref{lem:entropy-of-integral-operators}
\begin{eqnarray*}
\int_{0}^{\sigma^*_\mathbf{s}}\sqrt{\mathfrak{E}_{||\cdot||_2,\cQ_{\mathbf{s},\mathrm{x}}}(\delta)}\rd \delta&=&
\int_{0}^{\delta_0}\sqrt{\mathfrak{E}_{||\cdot||_2,\cQ_{\mathbf{s},\mathrm{x}}}(\delta)}\rd \delta+\int_{\delta_0}^{\sigma^*_\mathbf{s}}\sqrt{\mathfrak{E}_{||\cdot||_2,\cQ_{\mathbf{s},\mathrm{x}}}(\delta)}\rd \delta
\\
&\leq& A(\omega_2)\big(1-[2\omega_2]^{-1}\big)\delta_0^{1-\frac{1}{2\omega_2}}+ A(\omega_1)\big([2\omega_1]^{-1}-1\big)\delta_0^{1-\frac{1}{2\omega_1}}
\\*[2mm]
&=&\sqrt{\lambda^*\big(\omega_2,\mu\big)}\big(1-[2\omega_2]^{-1}\big)R_\mu^{\frac{1}{2\omega_2}} +
\sqrt{\lambda^*\big(\omega_1,\mu\big)}\big([2\omega_1]^{-1}-1\big)R_\mu^{\frac{1}{2\omega_1}}.
\end{eqnarray*}
It yields together with (\ref{eq142:proof-th:deviation-l_p-gauss})
$
 D_{\cQ_{\mathbf{s},\mathrm{x}},\rho}\leq C_\mu,
$
where, recall,
$$
C_\mu=4\sqrt{2}\|\cK\|_2^{d-1}\inf_{\{\omega_1,\omega_2\}\in\Omega}\bigg[\sqrt{\lambda^*\big(\omega_2,\mu\big)}\big(1-[2\omega_2]^{-1}\big)R_\mu^{\frac{1}{2\omega_2}} +
\sqrt{\lambda^*\big(\omega_1,\mu\big)}\big([2\omega_1]^{-1}-1\big)R_\mu^{\frac{1}{2\omega_1}}\bigg].
$$
Applying the assertion II of Lemma \ref{lem:talagrand-lifshits} we get
\begin{eqnarray}
\label{eq143:proof-th:deviation-l_p-gauss}
&& \bE\varsigma_{\mathbf{s}}(\mathrm{x})=
\bE\big\{\sup_{Q\in\cQ_{\mathbf{s},\mathrm{x}}}\varsigma_{\mathbf{s}}\big(Q,\mathrm{x}\big)\big\}\leq C_\mu.
\end{eqnarray}
We  obtain from (\ref{eq605:proof-th:deviation-l_p-gauss}) that
\begin{eqnarray}
\label{eq1430:proof-th:deviation-l_p-gauss}
\sigma_{\varsigma}:=\sup_{\mathbf{s}\in\cS_d}\sup_{Q\in\cQ_{\mathbf{s},\mathrm{x}}}\sqrt{\bE\varsigma^2_{\mathbf{s}}\big(Q,\mathrm{x}\big)}
=:\sup_{\mathbf{s}\in\cS_d}\sigma_{\mathbf{s}}
\leq  \|\cK\|_2^{d-1}\|\cK\|_{\frac{2\mu}{3\mu-2}}\mh^{\frac{1-\tau}{r}}.
\end{eqnarray}
Applying the assertion I of Lemma \ref{lem:talagrand-lifshits} we obtain in view
of (\ref{eq143:proof-th:deviation-l_p-gauss}) for any $z>0$
\begin{eqnarray}
\label{eq1431:proof-th:deviation-l_p-gauss}
&& \bP\Big\{\varsigma_{\mathbf{s}}(\mathrm{x})\geq C_\mu+z\Big\} \leq e^{-\frac{z^{2}}{2\sigma_\mathbf{s}^2}}\leq e^{-\frac{z^{2}}{2\sigma_\varsigma^2}}.
\end{eqnarray}
Set $U=C_\mu+\sqrt{2e^{r}}\|\cK\|_2^{d-1}\|\cK\|_{\frac{2\mu}{3\mu-2}}$ we obtain using (\ref{eq1431:proof-th:deviation-l_p-gauss})
\begin{eqnarray*}
\bE\varsigma(\mathrm{x})&\leq& U+\int_{0}^\infty
\bP\big\{\varsigma(\mathrm{x})\geq U+ y\big\}\rd y\leq
U + (S+1)^d\int_{0}^\infty e^{-\frac{(U-C_\mu+y)^{2}}{2\sigma_\varsigma^2}}  \rd y
\nonumber\\*[2mm]
 &\leq&
U+\sqrt{8\pi}\|\cK\|_2^{d-1}\|\cK\|_{\frac{2\mu}{3\mu-2}}(S+1)^d\exp{\Big\{-e^{r}\mh^{\frac{2(\tau-1)}{r}}\Big\}}.
\end{eqnarray*}
Taking into account that $(S+1)^d\leq \big[4\ln(\cA)\big]^d$ in view of the definition of $S$ and that
$$
\inf_{r>0}e^{r}\mh^{\frac{2(\tau-1)}{r}}=e^{2\sqrt{2(1-\tau)|\ln(\mh)|}},
$$
we obtain
\begin{eqnarray}
\label{eq1432:proof-th:deviation-l_p-gauss}
\bE\varsigma(\mathrm{x})&\leq& U+
\sqrt{8\pi}\|\cK\|_2^{d-1}\|\cK\|_{\frac{2\mu}{3\mu-2}}\big[4\ln(\cA)\big]^de^{e^{-2\sqrt{2(1-\tau)|\ln(\mh)|}}}.
\nonumber\\*[2mm]
 &\leq& C_\mu+
4^d\big(\sqrt{2e^{r}}+\sqrt{8\pi}\big)\|\cK\|_2^{d-1}\|\cK\|_{\frac{2\mu}{3\mu-2}}=\widetilde{C}_\mu.
\end{eqnarray}
The last inequality follows from the relation (\ref{eq0:relations-between-parameters}) and the definition of $U$.

\smallskip

$1^0\mathbf{b4}.\;$
Applying the assertion I of Lemma \ref{lem:talagrand-lifshits} we obtain
in view of (\ref{eq1430:proof-th:deviation-l_p-gauss})  for any $z>0$
\begin{eqnarray*}
&& \bP\Big\{\varsigma(\mathrm{x})\geq  \widetilde{C}_\mu+z\Big\} \leq e^{-\frac{z^{2}}{2\sigma_\varsigma^2}}.
\end{eqnarray*}
It yields together with (\ref{eq1430:proof-th:deviation-l_p-gauss})
\begin{eqnarray}
\label{eq603:proof-th:deviation-l_p-gauss}
\bE\Big(\varsigma^{\frac{r}{1-\tau}}(\mathrm{x})\Big)&=&\widetilde{C}_\mu^{\frac{r}{1-\tau}}+
 \frac{r}{1-\tau}\int_{0}^\infty \big(z+\widetilde{C}_\mu\big)^{\frac{r+\tau-1}{1-\tau}}
\bP\big\{\varsigma(\mathrm{x})\geq z+\widetilde{C}_\mu\big\}\rd z
\nonumber\\*[2mm]
&\leq& \widetilde{C}_\mu^{\frac{r}{1-\tau}}+\widehat{C}^{\frac{r}{1-\tau}}_\mu.
\end{eqnarray}
Here recall
$$
\widehat{C}_\mu=\bigg[\frac{r}{1-\tau}\int_{0}^\infty \big(u+\widetilde{C}_\mu\big)^{\frac{r+\tau-1}{1-\tau}}
 \exp{\Big\{-u^{2}\Big[2|\cK\|_2^{d-1}\|\cK\|_{\frac{2\mu}{3\mu-2}}\Big]^{-1}\Big\}}\rd u\bigg]^{\frac{1-\tau}{r}}.
$$
Similarly we deduce from (\ref{eq1430:proof-th:deviation-l_p-gauss}) and (\ref{eq1431:proof-th:deviation-l_p-gauss})
\begin{eqnarray}
\label{eq60400:proof-th:deviation-l_p-gauss}
\bE\Big(\varsigma_\mathbf{s}^{\frac{r}{1-\tau}}(\mathrm{x})\Big)
\leq C_\mu^{\frac{r}{1-\tau}}+\widehat{C}^{\frac{r}{1-\tau}}_\mu\leq \tilde{C}_\mu^{\frac{r}{1-\tau}}+\widehat{C}^{\frac{r}{1-\tau}}_\mu
,\quad\forall \mathbf{s}\in\bN^d.
\end{eqnarray}
Noting that the bounds in (\ref{eq603:proof-th:deviation-l_p-gauss}) and (\ref{eq60400:proof-th:deviation-l_p-gauss}) are independent of
$\mathrm{x}$ and $\mathbf{s}$ we get  in view of (\ref{eq602:proof-th:deviation-l_p-gauss})
\begin{eqnarray*}
\bE\big\{\zeta(r)\big\}\leq [1\vee (2b)^{d-1}\big]
\big[\cL^{\frac{1}{r}}+\cL^{\frac{\tau}{r}}\big(1-e^{-\frac{\tau p}{4}}\big)^{\frac{\tau-1}{r}}\big]
\big[\widetilde{C}_\mu+\widehat{C}_\mu\big].
\end{eqnarray*}
This proves (\ref{eq600:proof-th:deviation-l_p-gauss}) with
$T=[1\vee (2b)^{d-1}\big]
\big[\cL^{\frac{1}{r}}+\cL^{\frac{\tau}{r}}\big(1-e^{-\frac{\tau p}{4}}\big)^{\frac{\tau-1}{r}}\big]
\big[\widetilde{C}_\mu+\widehat{C}_\mu\big]$.

\smallskip

$1^0\mathbf{c}.\;$
Remembering that $C_2(r)=T
+e^r\sqrt{2(1+q)}\big(r\sqrt{e}\big)^d\|\cK\|_{\frac{2r}{r+2}}^d$ we obtain, applying the assertion I of Lemma \ref{lem:talagrand-lifshits}
available in view of (\ref{eq600:proof-th:deviation-l_p-gauss}) and (\ref{eq145:proof-th:deviation-l_p-gauss})
$$
\bP\Big\{\zeta(r)\geq C_2(r)+z\Big\}\leq e^{-e^{r}}e^{-qe^r\mh^{\frac{2d}{r}}}
\exp{\Big\{-\big[2\big(r\sqrt{e}\big)^d\|\cK\|_{\frac{2r}{r+2}}^d\big]^{-1}z^2\Big\}},\quad\forall z\geq 0.
$$
Taking into account that $e^r\mh^{\frac{2d}{r}}\leq e^{2\sqrt{2d|\ln(\mh)|}}$ for any $r>0$ we come to
 (\ref{eq6:proof-th:deviation-l_p-gauss}).

\smallskip

$2^0.\;$ We deduce from  (\ref{eq0002002:proof-th:deviation-l_p-gauss}) and (\ref{eq6:proof-th:deviation-l_p-gauss}) that
\begin{eqnarray*}
\label{eq146:proof-th:deviation-l_p-gauss}
\bE\big[\zeta(r)-C_2(r)\big]^q_+\leq \sqrt{(\pi/2)}e^{-e^{r}}\big[\big(r\sqrt{e}\big)^d
\|\cK\|_{\frac{2r}{r+2}}^d\big]^{\frac{q}{2}}e^{-qe^{2\sqrt{2d|\ln(\mh)|}}}\boldsymbol{\gamma}_{q+1},
\end{eqnarray*}
where recall $\boldsymbol{\gamma}_{q+1}$ is the $(q+1)$-th moment of the standard normal distribution.
This yields together with (\ref{eq00020022:proof-th:deviation-l_p-gauss})
\begin{eqnarray*}
\label{eq0002222:proof-th:deviation-l_p-gauss}
&&\bE\bigg(\sup_{\vec{h}\in\mathrm{H}}
\Big[\big\|\xi_{\vec{h}}\big\|_p-\psi_r(\vec{h})\Big]_+\bigg)^q
\leq  \bigg[C_4\cA e^{-e^{2\sqrt{2d|\ln(\mh)|}}}\bigg]^q,
\end{eqnarray*}
and the assertion of the theorem follows.

\epr

\subsection{Proof of Theorem \ref{th:deviation-l_p-gauss-isotropic}}

\subsubsection{Constants}
\label{sec:constants3}
Let $\mathbf{c}(d)$ be the constant appearing in   $(2,2)$-strong maximal inequality, see \cite{folland}. Set
\begin{eqnarray*}
\sigma_*&=&\sqrt{ 2^{d+1}a^d\|K\|_\infty\|K\|_1\mathbf{c}(d)}(2b)^{\frac{d(p-1)}{p}};
\\
C_5&=&\Big[\sqrt{8\pi}\sigma_*^{q-1}\boldsymbol{\gamma}_{q+1}\Big]^{\frac{1}{q}}\sum_{r=d+1}^\infty\sum_{l=1}^\infty e^{-2^{l}e^r}.
\end{eqnarray*}
For any $r\in\bN^*, r>d$ put $\gamma_r=\frac{d}{2}+\frac{d}{2pr}$ and  let  $\mathfrak{D}$ denote the unit disc in $\bR^d$.
Set
\begin{eqnarray*}
 T(r)&=&\big[\sigma_*/2\big]\vee\Big[(d/2+1)^dT^*(r)+\|\cK\|^d_1(2b)^{1/p}\Big];
\\
T^*(r)&=&2^{-d+1}\bigg[L(a+2)^{d}\int\mathfrak{z}^{-d-\gamma_r+\lfloor\gamma_r\rfloor+1}\mathrm{1}_{\mathfrak{D}}(\mathfrak{z})\rd \mathfrak{z}+
C(K)\int\mathfrak{z}^{-d-\gamma_r+\lfloor\gamma_r\rfloor}\mathrm{1}_{\overline{\mathfrak{D}}}(\mathfrak{z})\rd \mathfrak{z}\bigg],
\end{eqnarray*}
where $C(K)=\sup_{|\mathbf{n}|=\lfloor d/2\rfloor}\big\|D^{\mathbf{n}}K\big\|_1$.  Note that $\gamma_r\neq\lfloor\gamma_r\rfloor$ and, therefore, both integrals in the definition
of $T^*(r)$ are finite.

\smallskip

Let $\lambda_d^{*}(r)=\lambda_d\big(\gamma_r,1,1,[-a-b,a+b]^d\big)$,
where the quantity $\lambda_k(\cdot,\cdot,\cdot,\cdot), k\in\bN^*,$ is defined in Lemma \ref{lem:birman-solomyak}. Set finally
$$C_2^*(r)=8\sqrt{2\lambda_d^{*}(r)}\big[T(r)\big]^{d/2\gamma_r}
\big(\sigma_*/2\big)^{\frac{1}{2pr}}+4\sqrt{qe^r}\sigma_*.
$$

\subsubsection{Auxiliary lemma}

For any $l\in\bN^*$ and any $r\in\bN^*$ satisfying  $r>d$ put
$$
\mathrm{H}_{l,r}=\bigg\{\vec{h}\in\mathrm{H}:\;\; 2^{l-1}\mh^{-d}\leq\Big\|h^{-\frac{d}{2}}\Big\|_{p+\frac{1}{r}}< 2^{l}\mh^{-d} \bigg\}
$$
and introduce
$$
\mathfrak{Q}_{l,r}=\bigg\{Q:\bR^d\to\bR:\;\;Q(\cdot)=\int_{(-b,b)^d}K_{\vec{h}}(\cdot-x)\vartheta(x)\nu_d(\rd x),\;
\vartheta\in \bB_{q,d},\;\vec{h}\in\mathrm{H}_{l,r}\bigg\},
$$
where $\bB_{q,d}=\{\vartheta:(-b,b)^d\to\bR: \;\; \|\vartheta\|_q\leq 1\}$ and $1/q=1-1/p$.

\begin{lemma}
\label{lem:entropy-int-operators-on-Rd}
For any $r,l\in\bN^*$, $r>d$ and any $\delta\in \Big(0, T(r)\big(2^{l}\mh^{-d}\big)^{\frac{2\gamma_r}{d}}\Big]$ one has
$$
\mE_{\mathfrak{Q}_{l,r},\|\cdot\|_2}(\delta)\leq \lambda_d^{*}(r)\big[T(r)\big]^{d/\gamma_r}\big(2^{l}\mh^{-d}\big)^2\delta^{-d/\gamma_r}.
$$

\end{lemma}

\subsubsection{Proof of Theorem \ref{th:deviation-l_p-gauss-isotropic}} Set
$$
\psi^*_{r}(h)=C^*_2(r)\big\| h^{-\frac{d}{2}}\big\|_{p+\frac{1}{r}},\;r\in\bN^*,\; r>d.
$$
We have
$$
\bE\bigg\{\sup_{\vec{h}\in\mathrm{H}}
\Big[\big\|\xi_{\vec{h}}\big\|_p-\inf_{r\in\bN^*,r>d}\psi^*_r(\vec{h})\Big]_+\bigg\}^{q}
\leq \sum_{r=d+1}^\infty\bE\bigg\{\sup_{\vec{h}\in\mathrm{H}}
\Big[\big\|\xi_{\vec{h}}\big\|_p-\psi^*_r(\vec{h})\Big]_+\bigg\}^{q}.
$$
Moreover, since $\mathrm{H}=\cup_{l\geq 1}\mathrm{H}_{l,r}$ for any $r\in\bN^*$, one has
$$
\bigg\{\sup_{\vec{h}\in\mathrm{H}}
\Big[\big\|\xi_{\vec{h}}\big\|_p-\psi^*_r(\vec{h})\Big]_+\bigg\}^q\leq\sum_{l=1}^\infty
\bigg(\sup_{\vec{h}\in\mathrm{H}_{l,r}}
\big\|\xi_{\vec{h}}\big\|_p-C_2^*(r)2^{l-1}\mh^{-d}\bigg)^q_+.
$$
Thus,
\begin{eqnarray}
\label{eq1new:proof-th:deviation-l_p-gauss}
&&\bE\bigg\{\sup_{\vec{h}\in\mathrm{H}}
\Big[\big\|\xi_{\vec{h}}\big\|_p-\inf_{r\in\bN^*,r>d}\psi^*_r(\vec{h})\Big]_+\bigg\}^{q}
\leq \sum_{r=d+1}^\infty\sum_{l=1}^\infty\bE\bigg(\sup_{\vec{h}\in\mathrm{H}_{l,r}}
\big\|\xi_{\vec{h}}\big\|_p-C_2^*(r)2^{l-1}\mh^{-d}\bigg)_+^{q}.
\end{eqnarray}
We will proceed similarly to the proof of Theorem \ref{th:deviation-l_p-gauss}. First using   duality arguments we can assert that
$$
\sup_{\vec{h}\in\mathrm{H}_{l,r}}
\big\|\xi_{\vec{h}}\big\|_p=
\sup_{\vec{h}\in\mathrm{H}_{l,r}}\sup_{\vartheta\in\bB_{q,d}}\int_{(-b,b)^d} \xi_{\vec{h}}(x)\vartheta(x)\nu_d(\rd x).
$$
Noting that
$$
\int_{(-b,b)^d} \xi_{\vec{h}}(x)\vartheta(x)\nu_d(\rd x)=\int
\bigg[\int_{(-b,b)^d} h^{-d}(x) K\bigg(\frac{t-x}{h(x)}\bigg)\vartheta(x)\nu_d(\rd x)\bigg]W(\rd t)
$$
we have
$$
\sup_{\vec{h}\in\mathrm{H}_{l,r}}
\big\|\xi_{\vec{h}}\big\|_p=\sup_{Q\in\mathfrak{Q}_{l,r}}\int Q(t) W(\rd t)=:\sup_{Q\in\mathfrak{Q}_{l,r}}\zeta(Q).
$$
Thus, we get from (\ref{eq1new:proof-th:deviation-l_p-gauss})
\begin{eqnarray}
\label{eq100-new:proof-th:deviation-l_p-gauss}
\bE\bigg\{\sup_{\vec{h}\in\mathrm{H}}
\Big[\big\|\xi_{\vec{h}}\big\|_p-\inf_{r\in\bN^*,r>d}\psi^*_r(\vec{h})\Big]_+\bigg\}^{q}
\leq \sum_{r=d+1}^\infty\sum_{l=1}^\infty\bE\bigg(\sup_{Q\in\mathfrak{Q}_{l,r}}
\zeta(Q)-C_2^*(r)2^{l-1}\mh^{-d}\bigg)_+^{q}.
\end{eqnarray}
Obviously
$\zeta(\cdot)$
is centered gaussian random function on $\mathfrak{Q}_{l,r}$
and our goal is to apply to it  the assertion I  of Lemma \ref{lem:talagrand-lifshits}. To do this we have to show that
\begin{eqnarray}
\label{eq2new:proof-th:deviation-l_p-gauss}
\bE\Big\{\sup_{Q\in\mathfrak{Q}_{l,r}}\zeta(Q)\Big\}\leq U_{l,r}
\end{eqnarray}
 for some $0<U_{l,r}<\infty$ and to compute
\begin{eqnarray}
\label{eq3-new:proof-th:deviation-l_p-gauss}
\sigma_{l,r}^{2}:=\sup_{Q\in\mathfrak{Q}_{l,r} }\int Q^2(t)\nu_d(\rd t).
\end{eqnarray}
It is important to realize   that  this programm, being the same as in the proof of Theorem \ref{th:deviation-l_p-gauss},
requires  completely different arguments. It is related to the fact that we consider the random field $\xi_{\vec{h}}$ itself and not its "normalized" version
$\sqrt{V_{\vec{h}}}\xi_{\vec{h}}$.

\smallskip

$1^0.\;$ We start with bounding the quantity $\sigma_{l,r}$.
Putting for any $x,y\in (-b,b)^d$
$$
R(x,y)=\int K\bigg(\frac{t-x}{h(x)}\bigg)K\bigg(\frac{t-y}{h(y)}\bigg)\nu_d(\rd t),
$$
we obtain for any $Q\in\mathfrak{Q}_{l,r}$
\begin{eqnarray*}
\int Q^2(t)\nu_d(\rd t)&=&\int\bigg[\int_{(-b,b)^d} h^{-d}(x) K\bigg(\frac{t-x}{h(x)}\bigg)\vartheta(x)\nu_d(\rd x)\bigg]^2\nu_d(\rd t)
\\
&=&\int_{(-b,b)^d}\int_{(-b,b)^d}h^{-d}(x)h^{-d}(y)\vartheta(x)\vartheta(y)R(x,y)\nu_d(\rd x)\nu_d(\rd y).
\end{eqnarray*}
Taking into account that $\text{supp}(K)\subseteq [-a,a]^d$ in view of Assumption \ref{ass:kernel-new} we get
\begin{eqnarray*}
&&|R(x,y)|\leq\big[h(x)\wedge h(y)\big] \|K\|_\infty\|K\|_1\mathrm{1}_{[-2a,2a]}\bigg(\frac{x-y}{h(x)\vee h(y)}\bigg).
\end{eqnarray*}
Hence, putting $\Upsilon=\|K\|_\infty\|K\|_1$, we obtain
\begin{eqnarray*}
&&\int Q^2(t)\nu_d(\rd t)
\\
&&\leq \Upsilon\int_{(-b,b)^d}\int_{(-b,b)^d}\big|\vartheta(x)\vartheta(y)\big|\big[h(x)\vee h(y)\big]^{-d}\mathrm{1}_{[-2a,2a]^d}
\bigg(\frac{x-y}{h(x)\vee h(y)}\bigg)\nu_d(\rd x)\nu_d(\rd y).
\end{eqnarray*}
It remains to note
\begin{eqnarray*}
&&\big[h(x)\vee h(y)\big]^{-d}\mathrm{1}_{[-2a,2a]^d}\bigg(\frac{x-y}{h(x)\vee h(y)}\bigg)
\\
&&\leq h^{-d}(x)\mathrm{1}_{[-2a,2a]^d}\bigg(\frac{x-y}{h(x)}\bigg)+
 h^{-d}(y)\mathrm{1}_{[-2a,2a]^d}\bigg(\frac{x-y}{ h(y)}\bigg)
\end{eqnarray*}
and, therefore,
\begin{eqnarray*}
\int Q^2(t)\nu_d(\rd t)&\leq& 2\Upsilon\int_{(-b,b)^d}|\vartheta(v)|\bigg[\int_{(-b,b)^d}
h^{-d}(v)\mathrm{1}_{[-2a,2a]^d}\bigg(\frac{u-v}{h(v)}\bigg)
|\vartheta(u)|\nu_d(\rd u)\bigg]\nu_d(\rd v)
\\
&\leq& 2^{d+1}a^d\Upsilon\int|\vartheta^*(v)|\sup_{\lambda>0}(2\lambda)^{-d}\bigg[\int_{\bR^d}
\mathrm{1}_{[-\lambda,\lambda]^d}\bigg(\frac{u-v}{\lambda}\bigg)
|\vartheta^*(u)|\nu_d(\rd u)\bigg]\nu_d(\rd v)
\\
&\leq& 2^{d+1}a^d\Upsilon\int|\vartheta^*(v)|M[|\vartheta^*|](v)\nu_d(\rd v).
\end{eqnarray*}
Here we have put $\vartheta^*(\cdot)=\mathrm{1}_{(-b,b)^d}(\cdot)\vartheta(\cdot)$ and $M[|\vartheta^*|]$ denotes
the Hardy-Littlewood maximal operator applied to the function $|\vartheta^*|$.

In view of  $(2,2)$-strong maximal inequality, \cite{folland},  there exists $\mathbf{c}(d)$ such that
$$
\int_{\bR^d}\big\{M[|\vartheta^*|](v)\big\}^2\nu_d(\rd v)\leq \mathbf{c}^2(d)\int_{\bR^d}|\vartheta^*(v)|^2\nu_d(\rd v).
$$
Using the latter bound we obtain  applying  Cauchy-Schwartz inequality
\begin{eqnarray*}
\bigg[\int Q^2(t)\nu_d(\rd t)\bigg]^{\frac{1}{2}}\leq \sqrt{\mathbf{c}(d)}\bigg[\int_{(-b,b)^d} |\vartheta(v)|^2\nu_d(\rd v)\bigg]^{\frac{1}{2}}\leq
\sqrt{ 2^{d+1}a^d\Upsilon\mathbf{c}(d)}(2b)^{\frac{d(p-1)}{p}}.
\end{eqnarray*}
To get the last inequality we applied the H\"older inequality and took into account that $\vartheta\in\bB_{q,d}$ and $q\geq 2$ since $p\leq 2$.

Noting that the right hand side of the obtained inequality is independent of $Q$ we get
\begin{eqnarray}
\label{eq4-new:proof-th:deviation-l_p-gauss}
\sigma_{l,r}\leq
\sqrt{ 2^{d+1}a^d\|K\|_\infty\|K\|_1\mathbf{c}(d)}(2b)^{\frac{d(p-1)}{p}}:=\sigma_*.
\end{eqnarray}
We would like to emphasize that the condition $p\leq 2$ is crucial in order to obtain the
bound presented in (\ref{eq4-new:proof-th:deviation-l_p-gauss}).

\smallskip

$2^0.\;$ Let us now establish  (\ref{eq3-new:proof-th:deviation-l_p-gauss}). The intrinsic semi-metric $\rho_\zeta$ of $\zeta(\cdot)$
is given by
$$
\rho_\zeta(Q_1,Q_2)=\|Q_1-Q_2\|_2, \quad Q_1,Q_2\in\mathfrak{Q}_{l,r}.
$$
Taking into account that
$\frac{d}{2\gamma_r}=\frac{2pr}{2pr+1}<1$ and applying the second assertion of Lemma \ref{lem:talagrand-lifshits}
 and Lemma \ref{lem:entropy-int-operators-on-Rd} we obtain in view of (\ref{eq4-new:proof-th:deviation-l_p-gauss})
\begin{eqnarray*}
D_{\mathfrak{Q}_{l,r},\rho_{\zeta}}&=& 4\sqrt{2\lambda_d^{*}(r)}\big[T(r)\big]^{d/2\gamma_r}\big(2^{l}\mh^{-d}\big)
\int_{0}^{\sigma_{l,r}/2}\delta^{-d/2\gamma_r}\rd \delta
\nonumber\\
&=&4\sqrt{2\lambda_d^{*}(r)}\big[T(r)\big]^{d/2\gamma_r}\big(2^{l}\mh^{-d}\big)(\sigma_{l,r}/2)^{\frac{1}{2pr}}
\nonumber\\
&\leq&4\sqrt{2\lambda_d^{*}(r)}\big[T(r)\big]^{d/2\gamma_r}
\big(\sigma_*/2\big)^{\frac{1}{2pr}}\big(2^{l}\mh^{-d}\big).
\nonumber
\end{eqnarray*}
We conclude that Dudley integral is finite and as it is proved in Lemma \ref{lem:entropy-int-operators-on-Rd} $\mathfrak{Q}_{l,r}$ is a subset
 is totally bounded space with respect to the intrinsic semi-metric of $\zeta(\cdot)$. It implies that $\zeta(\cdot)$ is almost surely continuous on  $\mathfrak{Q}_{l,r}$
 and, therefore, $\sup_{Q\in\mathfrak{Q}_{l,r}}\zeta(Q)$ is a random variable.

 Thus, in view of the second assertion of  Lemma \ref{lem:talagrand-lifshits}
 \begin{eqnarray}
\label{eq5-new:proof-th:deviation-l_p-gauss}
\bE\Big\{\sup_{Q\in\mathfrak{Q}_{l,r}}\zeta(Q)\Big\}&\leq&4\sqrt{2\lambda_d^{*}(r)}\big[T(r)\big]^{d/2\gamma_r}
\big(\sigma_*/2\big)^{\frac{1}{2pr}}\big(2^{l}\mh^{-d}\big)
\end{eqnarray}
and  (\ref{eq3-new:proof-th:deviation-l_p-gauss}) is proved with $U_{l,r}=4\sqrt{2\lambda_d^{*}(r)}\big[T(r)\big]^{d/2\gamma_r}
\big(\sigma_*/2\big)^{\frac{1}{2pr}}\big(2^{l}\mh^{-d}\big)$.

Moreover,  $\zeta(\cdot)$ is almost surely bounded  on  $\mathfrak{Q}_{l,r}$ and, therefore,   the first assertion of Lemma \ref{lem:talagrand-lifshits}
is applicable.

\smallskip

$3^0.\;$ Hence,  noting that
$C_2^*(r)=8\sqrt{2\lambda_d^{*}(r)}\big[T(r)\big]^{d/2\gamma_r}
\big(\sigma_*/2\big)^{\frac{1}{2pr}}+4\sqrt{qe^r}\sigma_*
$
we obtain
 \begin{eqnarray*}
\label{eq6-new:proof-th:deviation-l_p-gauss}
\bP\Big\{\sup_{Q\in\mathfrak{Q}_{l,r}}\zeta(Q)\geq 2^{l-1}\mh^{-d}C_2^*(r)+z\Big\}\leq \exp{\big\{-2^{l+1}q\mh^{-d} e^r\big\}}
e^{-\frac{z^{2}}{2\sigma_*^2}},\quad\forall z>0.
\end{eqnarray*}
It yields for any $q\geq 1$
\begin{eqnarray}
\label{eq7-new:proof-th:deviation-l_p-gauss}
\bE\bigg(\sup_{Q\in\mathfrak{Q}_{l,r}}
\zeta(Q)-C_2^*(r)2^{l-1}\mh^{-d}\bigg)_+^{q}
&=&q\int_{0}^\infty z^{q-1}\bP\Big\{\sup_{Q\in\mathfrak{Q}_{l,r}}\zeta(Q)\geq 2^{l-1}\mh^{-d}C_2^*(r)+z\Big\}
\nonumber\\*[2mm]
&\leq&\sqrt{8\pi}\sigma_*^{q-1}\boldsymbol{\gamma}_{q+1}\exp{\big\{-2^{l+1}q\mh^{-d} e^r\big\}}.
\end{eqnarray}
We deduce from (\ref{eq100-new:proof-th:deviation-l_p-gauss}) and (\ref{eq7-new:proof-th:deviation-l_p-gauss})
\begin{eqnarray*}
&&\bE\bigg\{\sup_{\vec{h}\in\mathrm{H}}
\Big[\big\|\xi_{\vec{h}}\big\|_p-\inf_{r\in\bN^*,r>d}\psi^*_r(\vec{h})\Big]\bigg\}^{q}_+
\leq \Big(C_5e^{\mh^{-d}}\Big)^{q},
\end{eqnarray*}
where, recall, $C_5=\Big[\sqrt{8\pi}\sigma_*^{q-1}\boldsymbol{\gamma}_{q+1}\Big]^{\frac{1}{q}}\sum_{r=d+1}^\infty\sum_{l=1}^\infty e^{-2^{l}e^r}$.

\epr

\smallskip

\section{Appendix}

\paragraph{Proof of Lemma \ref{lem:entropy-of-integral-operators}}
Recall that  $\mu^{-1}=q^{-1}+\tau r^{-1}$ and note that $2>\mu>1$ since $\tau<1$ and $r>2$. The proof of the  lemma is  mostly based on the   inclusion
\begin{eqnarray}
\label{eq17:proof-th:deviation-l_p-gauss}
\cQ_{\mathrm{x},\mathbf{s}}&\in&\bS_\mu^{\omega}\Big([-a-b,a+b],\tilde{R}_\mu\Big),\;\quad \forall \omega\in\big(1/\mu-1/2,1\big),
\end{eqnarray}
where  $\tilde{R}_\mu=\|\cK\|_1+2\left[5\big\{4L(a+1)\big\}^\mu+4\big\{2\|\cK\|_1\big\}^\mu(2-\mu)^{-1}\right]^{\frac{1}{\mu}}$.

First, we note that all functions from $\cQ_{\mathrm{x},\mathbf{s}}$ vanish outside the interval $\Delta=[-a-b,a+b]$ since $\cK$ is compactly supported on $[-a,a]$
and $\mh_{s_1}\leq \mh<1$.


Next, applying Young inequality we obtain for any  $Q\in\cQ_{\mathrm{x},\mathbf{s}}$
\begin{eqnarray}
\label{eq18:proof-th:deviation-l_p-gauss}
\big\|Q\big\|_{\bL_\mu(\Delta)}&=&\lambda^{-1}_{\vec{h},\mathbf{s}}(\mathrm{x})\bigg[\int_{\Delta}
\bigg|\int_{-b}^b\mh_{s_1}^{-1/2}
\cK\bigg(\frac{y-x_1}{\mh_{s_1}}\bigg)\ell(x_1)\mathrm{1}_{\Lambda_{\mathbf{s}}\big[\vec{h}\big]}(x_1,\mathrm{x})\nu_1(\rd x_1)\bigg|^{\mu}
\nu_1(\rd y)\bigg]^{\frac{1}{\mu}}
\nonumber\\
&\leq& \lambda^{-1}_{\vec{h},\mathbf{s}}(\mathrm{x})(\mh_{s_1})^{\frac{1}{2}}\|\cK\|_1\bigg[\int_{-b}^b\left|\ell(x_1)\right|^\mu
\mathrm{1}_{\Lambda_{\mathbf{s}}\big[\vec{h}\big]}(x_1,\mathrm{x})\nu_1(\rd x_1)\bigg]^{\frac{1}{\mu}}\leq(\mh_{s_1})^{\frac{1}{2}}\|\cK\|_1.
\end{eqnarray}
To get the last inequality we have used   (\ref{eq604:proof-th:deviation-l_p-gauss}).

Let $\omega\in \big(1/\mu-1/2,1\big)$ be fixed. Let us bound from above  the quantity
\begin{eqnarray*}
J_\mu:=\int_{\Delta}\int_{\Delta}\frac{|Q(y)-Q(z)|^{\mu}}{|y-z|^{1+\mu\omega}}\rd y\rd z.
\end{eqnarray*}

Putting $y=u+v$ and $z=u-v$ we obtain by changing of variables
$$
J_\mu\leq 2^{-\mu\omega}\int_{-\infty}^\infty |v|^{-1-\mu\omega}\bigg[\int_{-\infty}^\infty|Q_s(u+v)-Q_s(u-v)|^\mu\rd u\bigg]\rd v
$$
Note also that
\begin{eqnarray*}
&&|Q_s(u+v)-Q_s(u-v)|
\\*[2mm]
&&\leq \lambda^{-1}_{\vec{h},\mathbf{s}}(\mathrm{x}) \int_{-b}^b\mh_{s_1}^{-1/2}\bigg|\cK\bigg(\frac{u-x_1}{\mh_{s_1}}+\frac{v}{\mh_{s_1}}\bigg)
-\cK\bigg(\frac{u-x_1}{\mh_{s_1}}-\frac{v}{\mh_{s_1}}\bigg)\bigg||\ell(x_1)|
\mathrm{1}_{\Lambda_{\mathbf{s}}\big[\vec{h}\big]}(x_1,\mathrm{x})\nu_1(\rd x_1).
\end{eqnarray*}
Hence,
$$
J_\mu\leq 2^{-\mu\omega}\mh_s^{-\mu(\omega+1/2)}\lambda^{-\mu}_{\vec{h},\mathbf{s}}(\mathrm{x})\int_{-\infty}^\infty |w|^{-1-\mu\omega}G^{\mu}(w)\rd w,
$$
where  we have put for any $w\in\bR$
$$
G(w)=\bigg[\int_{-\infty}^\infty\bigg[\int_{-b}^b\bigg|\cK\bigg(\frac{u-x_1}{\mh_{s_1}}+w\bigg)
-\cK\bigg(\frac{u-x_1}{\mh_{s_1}}-w\bigg)\bigg||\ell(x_1)|\mathrm{1}_{\Lambda_{\mathbf{s}}\big[\vec{h}\big]}(x_1,\mathrm{x})\nu_1(\rd x_1)
\bigg]^\mu\rd u\bigg]^{\frac{1}{\mu}}.
$$
Applying Young inequality for any fixed $w$ and we obtain
\begin{eqnarray*}
G(w)&\leq& \mh_{s_1}\bigg[\int_{-\infty}^\infty\left|\cK\big(u+w\big)-\cK\big(u-w\big)\right|\rd u\bigg]
\bigg[\int_{-b}^b|\ell(x_1)|^\mu\mathrm{1}_{\Lambda_{\mathbf{s}}\big[\vec{h}\big]}(x_1,\mathrm{x})\nu_1(\rd x_1)
\bigg]^{\frac{1}{\mu}}
\\
&\leq&\mh_{s_1}\bigg[\int_{-\infty}^\infty\left|\cK\big(u+w\big)-\cK\big(u-w\big)\right|\rd u\bigg]\lambda_{\vec{h},\mathbf{s}}(\mathrm{x}).
\end{eqnarray*}
To get the last inequality we have used   (\ref{eq604:proof-th:deviation-l_p-gauss}).
Note  that
\begin{eqnarray*}
\int_{-\infty}^\infty\left|\cK\big(u+w\big)-\cK\big(u-w\big)\right|\rd u&\leq &2\|\cK\|_1,\quad\forall w\in\bR;
\\
\int_{-\infty}^\infty\left|\cK\big(u+w\big)-\cK\big(u-w\big)\right|\rd u&\leq& 4L(a+1)|w|,\quad \forall w\in[-1,1].
\end{eqnarray*}
To get the second inequality we have used Assumption \ref{ass:kernel} ($\mathbf{i}$). Thus, we get finally
\begin{eqnarray}
\label{eq21:proof-th:deviation-l_p-gauss}
&&J_\mu\leq 2^{-\mu\omega}\mh_{s_1}^{\mu(1/2-\omega)}\left[5\big\{4L(a+1)\big\}^\mu+4\big\{2\|\cK\|_1\big\}^\mu(2-\mu)^{-1}\right].
\end{eqnarray}
Here we have also used that $\mu<2$ and $\mu\omega>(2-\mu)(2\mu)^{-1}$.

Putting $\tilde{R}_\mu=\|\cK\|_1+\left[5\big\{2L(a+2)\big\}^\mu+4\big\{2\|\cK\|_1\big\}^\mu(2-\mu)^{-1}\right]^{\frac{1}{\mu}}$ we get from
(\ref{eq18:proof-th:deviation-l_p-gauss}) and (\ref{eq21:proof-th:deviation-l_p-gauss})
for any $\omega\in \big(1/\mu-1/2,1\big)$
\begin{eqnarray*}
&&\big\|Q\big\|_{\bL_\mu(\Delta)}+\bigg[\int_{\Delta}\int_{\Delta}\frac{|Q(y)-Q(z)|^\mu}{|y-z|^{1+\mu\omega}}\rd y\rd z\bigg]^{1/\mu}\leq
\tilde{R}_\mu \mh_{s_1}^{\frac{1}{2}-\omega}.
\end{eqnarray*}
Thus, the inclusion (\ref{eq17:proof-th:deviation-l_p-gauss}) is proved since $\tilde{R}_\mu\leq R_\mu$.
The assertion of the lemma follows from Lemma \ref{lem:birman-solomyak} with $k=1$ and its consequence (\ref{eq:cor-birman-solomjak}).

\epr

\subsection{Proof of Lemma \ref{lem:entropy-int-operators-on-Rd}}

Similarly to the proof of Lemma \ref{lem:entropy-of-integral-operators} the proof of the present lemma is  based on the  inclusion
\begin{eqnarray}
\label{eq1:proof-lemm:SS-general}
\mathfrak{Q}_{l,r}\subset\bS^{\gamma_r}_1\Big((-a-b,a+b)^{^{d}},R\Big), \quad R=T(r)
\big(2^{l}\mh^{-d}\big)^{\frac{2\gamma_r}{d}}.
\end{eqnarray}
Indeed, if (\ref{eq1:proof-lemm:SS-general}) holds then the required assertion  follows from the consequence (\ref{eq:cor-birman-solomjak}) of Lemma \ref{lem:birman-solomyak}.

 Thus, let us prove (\ref{eq1:proof-lemm:SS-general}).
First, we note that all functions from $\mathfrak{Q}_{l,r}$ vanish outside  the cube $\boldsymbol{\Delta}=[-a-b,a+b]^d$ since $K$
is compactly supported on $[-a,a]^d$
and $\mh<1$.

Next, for any $Q\in\mathfrak{Q}_{l,r}$ we obviously have
\begin{eqnarray}
\label{eq01:proof-lemm:SS-general}
\|Q\|_1:=\int_{\boldsymbol{\Delta}}|Q(t)|\nu_d(\rd t)\leq \|\cK\|^d_1\int_{(-b,b)^d}|\vartheta(x)|\nu_d(\rd x)\leq \|\cK\|^d_1(2b)^{1/p},
\end{eqnarray}
where the last inequality follows from the condition $\vartheta\in\bB_{q,d}$ and the H\"older inequality.

Taking into account that $\vec{h}(x)=\big(h(x),\ldots,h(x)\big)$ and that $\lfloor \gamma_r\rfloor=\lfloor d/2\rfloor$,
 we have for any $\mathbf{n}\in\bN^d$ satisfying
$|\mathbf{n}|=\lfloor \gamma_r\rfloor$ in view of Assumption \ref{ass:kernel-new}
$$
D^{\mathbf{n}}Q(t)=\int_{(-b,b)^d}\big[h(x)\big]^{-|\mathbf{n}|-d}\big[D^{\mathbf{n}}K\big]\bigg(\frac{t-x}{\vec{h}(x)}\bigg)\vartheta(x)\nu_d(\rd x).
$$
Moreover, putting $y=u+v$ and $z=u-v$ we obtain by changing of variables
$$
I_{\mathbf{n}}:=\int_{\boldsymbol{\Delta}}\int_{\boldsymbol{\Delta}}
\frac{\big|D^{\mathbf{n}}Q(y)-D^{\mathbf{n}}Q(z)\big|}{|y-z|^{d+\alpha}}
\rd y\rd z\leq 2^{-d-\alpha}\int_{\boldsymbol{\bR^d}}|v|^{-d-\alpha}T(v)\rd v.
$$
Here $\alpha=\gamma_r-\lfloor\gamma_r\rfloor$ and $T(v)=\int_{\bR^d}\big|D^{\mathbf{n}}Q(u+v)-D^{\mathbf{n}}Q(u-v)\big|\rd u$.

We get using Fubini theorem
\begin{eqnarray*}
I_\mathbf{n}&\leq&  2^{-d-\alpha}\int_{(-b,b)^d}\big[h(x)\big]^{-|\mathbf{n}|-d}|\vartheta(x)|\bigg\{\int|v|^{-d-\alpha}
\bigg[\int
\\
&&
\bigg|\big[D^{\mathbf{n}}K\big]\bigg(\frac{u+v-x}{h(x)}\bigg)-
\big[D^{\mathbf{n}}K\big]\bigg(\frac{u-v-x}{h(x)}\bigg)\bigg|\rd u\bigg]\rd v\bigg\}\nu_d(\rd x),
\end{eqnarray*}
By changing variables in inner integrals  $w=(u-x)/h(x)$ and  $\mathfrak{z}=v/h(x)$ we obtain
\begin{eqnarray}
\label{eq2:proof-lemm:SS-general}
I_\mathbf{n}&\leq&  T\int_{(-b,b)^d}\big[h(x)\big]^{-|\mathbf{n}|-\alpha}|\vartheta(x)|\nu_d(\rd x),
\end{eqnarray}
where $T= 2^{-d-\alpha}\int|\mathfrak{z}|^{-d-\alpha}\int\big|D^{\mathbf{n}}K(w+\mathfrak{z})-
D^{\mathbf{n}}K(w-\mathfrak{z})\big|\rd w\rd \mathfrak{z}$.

\smallskip

We obtain in view of Assumption \ref{ass:kernel-new} for any $|\mathbf{n}|\leq \lfloor d/2\rfloor+1$
\begin{eqnarray*}
\int\left|D^{\mathbf{n}}K(w+\mathfrak{z})-D^{\mathbf{n}}K(w-\mathfrak{z})\right|\rd w&\leq &2C(K),\quad\forall \mathfrak{z}\in\bR^d;
\\
\int\left|D^{\mathbf{n}}K(w+\mathfrak{z})-D^{\mathbf{n}}K(w-\mathfrak{z})\right|\rd w&\leq& 2L(a+2)^{d}|\mathfrak{z}|,
\quad \forall |\mathfrak{z}|\leq 1.
\end{eqnarray*}
 It yields (recall that $\mathfrak{D}$ denotes the unit disc in $\bR^d$),
$$
T\leq  2^{-d+1}\bigg[L(a+2)^{d}\int\mathfrak{z}^{-d-\alpha+1}\mathrm{1}_{\mathfrak{D}}(\mathfrak{z})\rd z+
C(K)\int\mathfrak{z}^{-d-\alpha}\mathrm{1}_{\overline{\mathfrak{D}}}(\mathfrak{z})\rd \mathfrak{z}\bigg]=T^*(r).
$$
 Thus,  we deduce from (\ref{eq2:proof-lemm:SS-general}) for any $\mathbf{n}$ satisfying $|\mathbf{n}|=\lfloor \gamma_r\rfloor$
 \begin{eqnarray}
\label{eq3:proof-lemm:SS-general}
I_\mathbf{n}&\leq&  T^*(r)\int_{(-b,b)^d}\big[h(x)\big]^{-\gamma_r}|\vartheta(x)|\nu_d(\rd x)
\leq  T^*(r)\bigg(\int_{(-b,b)^d}\big[h(x)\big]^{-p\gamma}\nu_d(\rd x)\bigg)^{\frac{1}{p}}
\nonumber\\
&=&
 T^*(r) \bigg(\Big\|h^{-\frac{d}{2}}\Big\|_{\frac{2p\gamma_r}{d}}\bigg)^{\frac{2\gamma_r}{d}}
= T^*(r) \bigg(\Big\|h^{-\frac{d}{2}}\Big\|_{p+\frac{1}{r}}\bigg)^{\frac{2\gamma_r}{d}}.
\end{eqnarray}
Here we have used H\"older inequality, the condition $\vartheta\in\bB_{q,d}$ and the definition of $\gamma_r$.

Taking into account that $\vec{h}\in\mathrm{H}_{l,r}$ we obtain from (\ref{eq3:proof-lemm:SS-general}) that
$$
\sum_{|\mathbf{n}|=\lfloor \gamma_r\rfloor}I_\mathbf{n}\leq  (d/2+1)^dT^*(r)\big(2^{l}\mh^{-d}\big)^{\frac{2\gamma_r}{d}}.
$$
It leads together with (\ref{eq01:proof-lemm:SS-general}) to the assertion of the lemma.
\epr

\subsection{Proof of Proposition \ref{prop:set-of-bandw}} Set
$$
B_{\vec{h}}(f,x)=\bigg|\int K_{\vec{h}}(t-x)f(t)\rd t-f(x)\bigg|,\;\; x\in\bR^d.
$$
We start the proof with several remarks.

1) Obviously  $\Lambda_{\mathbf{s}}\big[\vec{h}_f\big]\in\mB\big(\bR^d\big)$ for any $f\in\bN_d(\vec{\beta},\vec{r},\vec{L})$ and any multi-index $\mathbf{s}$ since $B_{\vec{h}}(f,\cdot)$ is measurable function. Moreover $\vec{h}_{f}(\cdot)$ takes its values in countable set that implies that $\vec{h}_{f}(\cdot)$ is measurable function.

2) The definition of the Nikolskii class implies that $\|f\|_{r_j}\leq L_j$ for any $j=1,\ldots d$. It yields, in view of the Young inequality
$$
\big\|B_{\vec{h}}(f,\cdot)\big\|_{r_j}\leq (1+\|K\|_1)L_j,\quad \forall j=1,\ldots d,
$$
and therefore,
$$
\nu_d\big(x\in (-b,b)^d:\; B_{\vec{h}}(f,x)=\infty\big)=0,\quad \forall \vec{h}\in\mH^{d}_\e.
$$
This, in its turn, implies that
\begin{equation}
\label{eq1:proof-prop:set-of-bandw}
\nu_d\Big(\cup_{j=1}^d\big\{x\in (-b,b)^d:\;\; h_j(f,x)=\infty\big\}\Big)=0.
\end{equation}

3) The following statement was proved in \cite{GL13}, Lemma 3: there exists a constant $\widetilde{C}$ completely determined by
$\vec{\beta}, d$ and the function $w$ such that
\begin{equation}
\label{eq:bias-norms-1}
B_{\vec{h}}(f,x)\leq\sum_{j=1}^dB_{\vec{h},j}(f,x),\;\;x\in\bR^d,\qquad  \big\|B_{\vec{h},j}(f, \cdot)\big\|_{r_j} \leq \tilde{C}L_j h_j^{\beta_j},\;\;\;\forall j=1,\ldots,d.
\end{equation}

{\it $1^{0}$. Proof of the first assertion.}
For any $\mathbf{s}\in\bN^*$   recall that $\vec{\mh}_\mathbf{s}=(\mh_{s_1},\ldots,\mh_{s_d})$ and
$V_{\mathbf{s}}=\prod_{j=1}^d\mh_{s_j}$. Denote by $\cS_d$ the set consisting of $\mathbf{s}=(s_1,\ldots, s_d)\in\bN^d$ satisfying
$s_j\geq S_\e(j)$  for any $j=1,\ldots, d.$ We will also use the following notation: for any $\mathbf{s}\in\cS_d$ let $\hat{\mathbf{s}}\in\bN^d$
be such that $\hat{\mathbf{s}}<\mathbf{s}$ and $|\mathbf{s}-\hat{\mathbf{s}}|=1$.

\smallskip

Putting  $\cX=\cap_{j=1}^d\{x\in (-b,b)^d:\;h_j(f,x)<\infty\}$ we have  in view of the definition $\vec{h}(f,\cdot)$ for any
$\mathbf{s}\in\cS_d$ such that $\mathbf{s}\neq \big(S_e(1),\ldots, S_e(d)\big)$.
\begin{eqnarray*}
\Lambda_{\mathbf{s}}\big[\vec{h}_f\big]\cap\cX&\subseteq&\left\{x\in (-b,b)^d:\;B_{\vec{\mh}_\mathbf{s}}(f,x)+\e V^{-\frac{1}{2}}_{\mathbf{s}}\leq
B_{\vec{\mh}_{\mathbf{\hat{s}}}}(f,x)+\e V^{-\frac{1}{2}}_{\mathbf{\hat{s}}}\right\}
\\
&\subseteq&\left\{x\in (-b,b)^d:\;
B_{\vec{\mh}_{\mathbf{\hat{s}}}}(f,x)\geq \e V^{-\frac{1}{2}}_{\mathbf{s}}(1-e^{-1/2}) \right\}
\\
&\subseteq&\bigcup_{j=1}^d\left\{x\in (-b,b)^d:\;
B_{\vec{\mh}_{\mathbf{\hat{s}}},j}(f,x)\geq \e V^{-\frac{1}{2}}_{\mathbf{s}}(1-e^{-1/2})d^{-1} \right\}
\end{eqnarray*}
The last inclusion follows from the first inequality in (\ref{eq:bias-norms-1}) and the definition of $\mathbf{\hat{s}}$.

We get from
(\ref{eq1:proof-prop:set-of-bandw}),  the second inequality in  (\ref{eq:bias-norms-1}) and Markov inequality
\begin{eqnarray*}
\nu_d\Big(\Lambda_{\mathbf{s}}\big[\vec{h}_f\big]\Big)=\nu_d\Big(\Lambda_{\mathbf{s}}\big[\vec{h}_f\big]\cap\cX\Big)&\leq&
\sum_{j=1}^d d^{r_j}V^{\frac{r_j}{2}}_{\mathbf{s}}\big[\e(1-e^{-1/2})\big]^{-r_j}\big\|B_{\vec{\mh}_{\mathbf{\hat{s}}},j}(f,\cdot)\big\|^{r_j}_{r_j}
\\
&\leq&\sum_{j=1}^d \kappa_j \big[\e^{-1}V^{\frac{1}{2}}_{\mathbf{s}}\mh_{s_j}^{\beta_j}\big]^{r_j},
\end{eqnarray*}
where we have put $\kappa_j=\big\{d(e^{\beta_j}-e^{\beta_j-1/2})\widetilde{C}L_j\big\}^{r_j}$ and used once again the definition of $\mathbf{\hat{s}}$.

Since $\nu_d\Big(\Lambda_{\mathbf{s}}\big[\vec{h}_f\big]\Big)=0$ for any $\mathbf{s}\notin \cS_d$ by the definition of $\vec{h}_f$ and
$\nu_d\Big(\Lambda_{\mathbf{s_0}}\big[\vec{h}_f\big]\Big)\leq (2b)^{d}$, $\mathbf{s_0}=\big(S_\e(1),\ldots,S_\e(d)\big)$,
we obtain for any $\tau\in (0,1)$
\begin{eqnarray*}
\sum_{\mathbf{s}\in\bN^d}\nu^{\tau}_d\Big(\Lambda_{\mathbf{s}}\big[\vec{h}_f\big]\Big)\leq \sum_{j=1}^d \kappa^\tau_j \sum_{\mathbf{s}\in\cS_d, \mathbf{s}\neq\mathbf{s_0}} \big[\e^{-1}V^{\frac{1}{2}}_{\mathbf{s}}\mh_{s_j}^{\beta_j}\big]^{\tau r_j}+(2b)^{d}.
\end{eqnarray*}
In view of  (\ref{eq:example-def-S(j)}) (the definition of $S_\e(j), j=1,\ldots,d$)  we get
\begin{eqnarray*}
V^{\frac{1}{2}}_{\mathbf{s}}&=&\Big[\mh^{d}e^{-\sum_{l=1}^dS_\e(l)}e^{\sum_{l=1}^d(S_\e(l)-s_l)}\Big]^{\frac{1}{2}}\leq\e^{\frac{1}{2\beta+1}}
e^{\frac{1}{2}\sum_{l=1}^d(S_\e(l)-s_l)};
\\
\mh_{s_j}^{\beta_j}&=&\mh e^{-\beta_j S_\e(j)}e^{\beta_j(S_\e(j)-s_j)}\leq \e^{\frac{2\beta}{2\beta+1}}e^{\beta_j(S_\e(j)-s_j)}\leq \e^{\frac{2\beta}{2\beta+1}}.
\end{eqnarray*}
It yields
$
\e^{-1}V^{\frac{1}{2}}_{\mathbf{s}}\mh_{s_j}^{\beta_j}\leq e^{\frac{1}{2}\sum_{k=1}^d(S_\e(k)-s_k)}
 $
and, therefore,
\begin{eqnarray*}
\sum_{\mathbf{s}\in\bN^d}\nu^{\tau}_d\Big(\Lambda_{\mathbf{s}}\big[\vec{h}_f\big]\Big)\leq \sum_{j=1}^d \kappa^\tau_j
\Big(1-e^{-\frac{\tau r_j}{2}}\Big)^{-d}+(2b)^{d}=:\cL.
\end{eqnarray*}
The first assertion is proved.

\par

{\it $2^{0}$. Proof of the second assertion.} The condition of the proposition allows us to assert that there exists $\mathfrak{p}>p$ such that
$ \upsilon(2+1/\beta)> \mathfrak{p}$. Putting $\phi_\e= e^{d/2}\e^{\frac{2\beta}{2\beta+1}}$ we obtain using the definition of $\vec{h}_f$
\begin{eqnarray*}
\Big\|V^{-\frac{1}{2}}_{\vec{h}_f}\Big\|_{\mathfrak{p}}^{\mathfrak{p}}&\leq& \e^{-\mathfrak{p}}\Big\|B_{\vec{h}_{f}}(f,\cdot)+\e V^{-\frac{1}{2}}_{\vec{h}_f}\Big\|_{\mathfrak{p}}^{\mathfrak{p}}=\e^{-\mathfrak{p}}\int_{(-b,b)^d}\inf_{\vec{h}\in\mH_\e}\bigg[B_{\vec{h}}(f,x)+\e V^{-\frac{1}{2}}_{\vec{h}}\bigg]^{\mathfrak{p}}\rd x
\\
&\leq& \big(2\phi_\e\e^{-1}\big)^{\mathfrak{p}}+
\sum_{k=0}^\infty \big(2e^{k+1}\phi_\e\e^{-1}\big)^{\mathfrak{p}}\nu_d\bigg(x:\;\inf_{\vec{h}\in\mH_\e}\bigg[B_{\vec{h}}(f,x)+\e V^{-\frac{1}{2}}_{\vec{h}}\bigg]\geq 2e^{k}\phi_\e\bigg)
\\
&\leq& \big(2\phi_\e\e^{-1}\big)^{\mathfrak{p}}+
\sum_{k=0}^\infty \big(2e^{k+1}\phi_\e\e^{-1}\big)^{\mathfrak{p}}\nu_d\bigg(x:\;B_{\vec{\mh}[k]}(f,x)+\e V^{-\frac{1}{2}}_{\vec{\mh}[k]}\geq 2e^{k}\phi_\e\bigg),
\\
\end{eqnarray*}
where we choose $\vec{\mh}[k]\in\mH_\e$ as follows. Let  $\vec{h}[k]=\big(h_1[k],\ldots,h_d[k]\big)$ be given by
$$
h_j[k]=\big(\phi_\e)^{1/\beta_j}
e^{k\big(\frac{1}{\beta_j}-\frac{\upsilon(2+1/\beta)}{\beta_jr_j}\big)},\quad j=1,\ldots, d,
$$
and define $\vec{\mh}[k]\in\mH_\e$ from the relation $e^{-1}\vec{h}[k]\leq\vec{\mh}[k]<\vec{h}[k]$.

First we note that
$$
h_j[k]\leq \big(\phi_\e)^{1/\beta_j}\leq \mh e^{-S_\e(j)+1},
$$
since $\vec{r}\in [1,p]^d$ and $p<\upsilon(2+1/\beta)$. This guarantees the existence of  $\vec{\mh}[k]$.
Next,
$$
\e V^{-\frac{1}{2}}_{\vec{\mh}[k]}\leq \e V^{-\frac{1}{2}}_{e^{-1}\vec{h}[k]}=e^{k+d/2}\e^{\frac{2\beta}{2\beta+1}}=e^k\phi_\e,
$$
and, therefore, using the latter bound, (\ref{eq:bias-norms-1}) and Markov inequality we obtain

\begin{eqnarray*}
\Big\|V^{-\frac{1}{2}}_{\vec{h}_f}\Big\|_{\mathfrak{p}}^{\mathfrak{p}}
&\leq& \big(2\phi_\e\e^{-1}\big)^{\mathfrak{p}}+
\sum_{k=0}^\infty \big(2e^{k+1}\phi_\e\e^{-1}\big)^{\mathfrak{p}}\nu_d\bigg(x:\;B_{\vec{\mh}_{\mathbf{s}[k]}}(f,x)\geq e^{k}\phi_\e\bigg)
\\
&\leq& \big(2\phi_\e\e^{-1}\big)^{\mathfrak{p}}+
\sum_{k=0}^\infty \big(2e^{k+1}\phi_\e\e^{-1}\big)^{\mathfrak{p}}\sum_{j=1}^d (e^{k}\phi_\e)^{-r_j} (\tilde{C}L_j)^{r_j} \big(\mh_{s_j}[k]\big)^{\beta_jr_j}
\\
&\leq& \big(2\phi_\e\e^{-1}\big)^{\mathfrak{p}}+
\sum_{k=0}^\infty \big(2e^{k+1}\phi_\e\e^{-1}\big)^{\mathfrak{p}}e^{-k\upsilon(2+1/\beta)}\sum_{j=1}^d  \big(\tilde{C}L_j\big)^{r_j} .
\\
&=& \e^{-\frac{\mathfrak{p}}{2\beta+1}}\bigg\{\big(2e^{d/2}\big)^{\mathfrak{p}}+ \big(2e^{d/2+1}\big)^{\mathfrak{p}}
\sum_{k=0}^\infty e^{-k[\upsilon(2+1/\beta)-\mathfrak{p}]}\sum_{j=1}^d  \big(\tilde{C}L_j\big)^{r_j}
\bigg\}.
\end{eqnarray*}
As we see the assumption of the proposition $ \upsilon(2+1/\beta)> p$ allowing us to choose $\mathfrak{p}>p$ and $ \upsilon(2+1/\beta)> \mathfrak{p}$
is crucial. The second assertion is proved.

\epr

\par\bigskip

\bibliographystyle{agsm}

\end{document}